\documentclass{amsproc}
\usepackage{latexsym,amsmath}
\usepackage{amssymb}
\usepackage{amscd}

\newtheorem{theorem}{Theorem}
\newtheorem{lemma}{Lemma}
\newtheorem{convention}{Convention}
\newtheorem{corollary}{Corollary}
\numberwithin{equation}{section}

\begin{document}

\date{}
\author{M.I.Belishev}
\address{St. Petersburg   Department   of   V.A. Steklov    Institute   of   Mathematics
of   the   Russian   Academy   of   Sciences, 7, Fontanka, 191023
St. Petersburg, Russia and Saint Petersburg State University,
St.Petersburg State University, 7/9 Universitetskaya nab., St.
Petersburg, 199034 Russia.} \email{belishev@pdmi.ras.ru}\

\author{V.S.Mikhailov}
\address{St.Petersburg   Department   of   V.A.Steklov    Institute   of   Mathematics
of   the   Russian   Academy   of   Sciences, 7, Fontanka, 191023
St. Petersburg, Russia} \email{ftvsm78@gmail.com}

\title{Inverse problem for one-dimensional dynamical Dirac system (BC-method)} \maketitle

 \begin{abstract}
A forward problem for the Dirac system is to find
$u=\begin{pmatrix}u_1(x,t)\\u_2(x,t)\end{pmatrix}$ obeying
$iu_t+\begin{pmatrix}0&1\\-1&0\end{pmatrix}u_x+\begin{pmatrix}p&q\\q&-p\end{pmatrix}u=0$
for $x>0,\,t>0$;\,\,$u(x,0)=\begin{pmatrix}0\\0\end{pmatrix}$ for
$x {\geqslant} 0 $, and $u_1(0,t)=f(t)$ for $t>0$, with the real
$p=p(x), q=q(x)$. An input--output map $R: u_1(0,\cdot)\mapsto
u_2(0,\cdot)$ is of the convolution form $Rf=if+r\ast f$, where
$r=r(t)$ is  a {\it response function}. By hyperbolicity of the
system, for any $T>0$, function $r\big|_{0 {\leqslant} t
{\leqslant} 2T}$ is determined by $p,q\big|_{0 {\leqslant} x
{\leqslant} T}$.

An inverse problem is: for an (arbitrary) fixed $T>0$, given
$r\big|_{0 {\leqslant} t {\leqslant} 2T}$ to recover $p,q\big|_{0
{\leqslant} x {\leqslant} T}$. The procedure that determines $p,q$
is proposed, and the characteristic solvability conditions on $r$
are provided.

Our approach is purely time-domain and is based on studying the
controllability properties of the Dirac system. In itself the
system is not controllable: the local completeness of states does
not hold, but its relevant extension gains controllability. It is
the fact, which enables one to apply the boundary control method
for solving the inverse problem.
 \end{abstract}

\setcounter{section}{-1} \setcounter{equation}{-1}

\section{Introduction}\label{sec Introduction}
\subsubsection*{Motivation}
The Dirac system under consideration is one of basic models of
one-dimen\-si\-o\-nal quantum mechanics. For applications of the
Dirac-type systems to optics, see  \cite{Rak} and the papers cited
therein. In its original form, it is an {\it evolution} system of
hyperbolic type. As such, it possesses a fundamental property: the
finiteness of the domain of influence (FDI) holds and corresponds
to physical fact that sharp signals propagate in the system with
finite velocity. However, the known results on the inverse
problems (IPs), beginning from the classical ones
\cite{GasLev,Fadd}, mainly concern to the spectral
(frequency-domain) setups \cite{G1,G2,AS1,AS3}. It is implied that
time-domain problems can be reduced to spectral ones by the
Fourier transform. This is true but under an important remark:
passing to Fourier images, one looses the local character of
dependence of solutions on coefficients of the system. Such a
locality follows from FDI and is an intrinsic feature of
hyperbolic problems. It leads to a local dependence of
coefficients on the observation data, and motivates the relevant
(so-called time-optimal) setups of IPs. In the mean time, the use
of such a locality is a specific feature and advantage of the {\it
boundary control method}, which is an approach to IPs based on
their relations to control theory \cite{BIP07}, \cite{1D-BC}
\footnote{for another approaches keeping the locality principle
see \cite{Blag,Nizh,Rak}}. The goal of the paper is to extend the
BC-method to the {\it time-domain inverse problem} for the Dirac
system.

\subsubsection*{Content and results}
{\bf Forward problem.}\,\,\,In its original form, a Dirac system
is associated with the initial boundary-value problem
\begin{align*}
& iu_t+Ju_x+Vu=0, && x>0, \quad 0<t<T \\
& u\big|_{t=0}=0, && x \geqslant 0 \\
& u_1\big|_{x=0}=f, && 0\leqslant t \leqslant T\,,
\end{align*}
where $i=\sqrt{-1}$\,; $J:=\begin{pmatrix}
0&1\\-1&0\end{pmatrix}$; $V=\begin{pmatrix} p&q\\q&-p
\end{pmatrix}$ is a matrix {\it potential}, $p=p(x)$ and $q=q(x)$ being
real-valued functions; $T>0$ is a final moment; $f=f(t)$ is a
complex-valued function ({\it boundary control});
$u=u^f(x,t)=\begin{pmatrix} u_1(x,t)\\
u_2(x,t) \end{pmatrix}$ is a solution ({\it state}).
Well-possedness of this problem and appropriate representation of
its solution are established by reducing to relevant integral
equations. Hyperbolicity of the problem implies $u^f\big|_{t<x}=0$
that corresponds to FDI. As a consequence, a {\it local extended
problem} of the form
\begin{align*}
 & iu_t+Ju_x+Vu=0, && 0<x<T, \quad 0<t<2T-x \\
 & u\big|_{t<x}=0\\
 & u_1\big|_{x=0}=f, && 0\leqslant t \leqslant 2T\,,
 \end{align*}
turns out to be well posed.
\smallskip

\noindent{\bf Dynamical system.}\,\,\,With the latter problem one
associates a dynamical system, which is denoted by $\alpha^T$ and
endowed with the proper control theory attributes: spaces and
operators. By FDI, its {\it reachable sets} ${{\mathcal
U}}^t_\alpha=\{u^f(\cdot,t)\,|\,\,f \in L_2([0,2T]; {\mathbb
C})\}\,\,(0 < t \leqslant T)$ are embedded to $L_2([0,t]; {\mathbb
C}^2)$. However, this embedding {\it is not dense}: $L_2([0,t];
{\mathbb C}^2)\ominus {{\mathcal U}}^t_\alpha\not=\{0\}$ holds for
any $t$. This fact is interpreted as a lack of the {\it local
boundary controllability}. By philosophy of the BC-method, such a
lack always complicates analysis and solving inverse problems, and
should be eliminated. By doing so, we introduce an auxiliary
system $\bar \alpha^T$ governed by the dual equation
$iv_t-Jv_x-Vv=0$ and compose a system $\beta^T$, which, with some
abuse of terms, may be specified as $\alpha^T \oplus \bar
\alpha^T$. System $\beta^T$ turns out to be controllable, and it
is the fact, which opens the way to apply the BC-method.
\smallskip

\noindent{\bf Inverse problem.}\,\,\,In the system $\alpha^T$, an
input-output correspondence is rea\-lized by a {\it response
operator} associated with the local extended  problem. It acts in
$L_2([0,2T]; {\mathbb C})$ by $R f:=u^f_2(0,\cdot)\big|_{0
{\leqslant} t {\leqslant} 2T}$ and is represented in the
convolution form $Rf=if+r\ast f$, where $r\big|_{0 {\leqslant} t
{\leqslant} 2T}$ is a {\it response function}. By FDI, the
response function is determined by the values of the potential on
$0{\leqslant} x {\leqslant} T$ only (doesn't depend on
$V\big|_{x>T}$). Therefore, the relevant to FDI setup of the
inverse problem is: for a fixed $T>0$, given $r\big|_{0
{\leqslant} t {\leqslant} 2T}$, to recover $V\big|_{0{\leqslant} x
{\leqslant} T}$.

To formulate the main result (Theorem \ref{Th Charact}), we accept
the following agreement on notation and keep it throughout the
paper.
 \begin{convention}\label{Conv 1}
All functions of time $t\geqslant 0$ are extended to $t<0$ by
zero.
 \end{convention}
Also, for a $z \in \mathbb C$, $\bar z:=\Re z -i \Im z$ is its
conjugate.

A $\mathbb C$-valued function $r\big|_{0 {\leqslant} t {\leqslant} 2T}$ determines
an operator $C$ acting in $L_2([0,2T]; {\mathbb C}^2)$ by the rule

\begin{equation*}\label{cal C in detail} \left(Ca\right)(t)\,=\,2 a(t)+\int_0^T
c(t,s) a(s)\,ds\,, \qquad 0\leqslant t \leqslant T
 \end{equation*}
for $a=\begin{pmatrix}a_1(t)\\a_2(t)\end{pmatrix}$, with the
matrix kernel elements
 \begin{align*}
\notag & c_{11}(t,s)(t)=-i\,[r(t-s)-\bar r(s-t)]\,, \quad
c_{12}(t,s)=-i\,\bar r(2T-t-s)\,,\\
& c_{21}(t,s)= i\,r(2T-t-s)\,,\quad c_{22}(t,s)(t)=i\,[\bar
r(t-s)-r(s-t)]\,.
 \end{align*}
We prove that $r\in C^1([0,2T]; {\mathbb C})$ is the response
function of a certain system $\alpha^T$ with a {\it $C^1$-smooth
real zero trace} potential $V$ if and only if operator $C$ is a
positive definite isomorphism \footnote{An {\it isomorphism} is a
surjective, bounded, and boundedly invertible linear operator. A
{\it positive definiteness} means that $(Ca,a)\geqslant
\theta\|a\|^2$ holds for all $a$ with a constant $\theta>0$.} in
$L_2([0,2T]; {\mathbb C}^2)$. Also, an efficient procedure which
recovers $V$ from $r$ is provided.

\subsubsection*{Comments}
$\bullet$\quad Solving the IP, we invoke no spectral devices (like
Fourier transform) and don't leave the finite space interval $0
{\leqslant} x {\leqslant} T$. Such an option is provided by the
FDI-principle and is one of the advantages of the BC-method.
\smallskip

\noindent $\bullet$\quad As a solvability condition for IP, a
positive definiteness of the operator $C$ is quite natural and
expectable. The novelty is its concrete form for the Dirac system.
Remarkably, all the rich properties of the system turn out to be
encoded in the specific structure of the matrix elements $c_{kl}$
and can be extracted from a single function $r$.
\smallskip

\noindent $\bullet$\quad We deal with $C^1$-smooth potentials and
inverse data just because this enables one to use classical
solutions that is convenient from the technical viewpoint.
Undoubtedly, our results can be extended to a wider class of
potentials.
\smallskip

\noindent $\bullet$\quad The given paper inscribes one more
dynamical IP in the framework of the BC-method. Its basic theses
are the following:

$\star$\,\,\,To solve a dynamical IP, one has to make clear the
character of controllability of the corresponding dynamical
system. By the very general principle of system theory, the better
is the controllability, the richer is information about the
system, which can be extracted from the inverse data
\cite{KFA,BIP07,1D-BC}.

$\star$\,\,\,The key object of our approach is the connecting
operator $C$ of the system. Roughly speaking, to solve the IP by
the BC-method is to express $C$ via the inverse data and then find
$C^{-1}$. In the mean time, to invert $C$ is to solve a family of
{\it linear} integral equations, which are a relevant analog of
the classical Gelfand-Levitan-Krein-Marchenko equations
\cite{BM_1}.

$\star$\,\,\,From the functional analysis viewpoint, the
background of the above-mentioned scheme is a {\it triangular
factorization} of the form $C=W^*W$. Such a factorization is
realized via constructing a family of projectors ${\mathcal
P}^\xi$ orthogonal in the relevant Hilbert metric determined by
$C$.
\smallskip

\noindent Such a scheme of solving IPs is reproduced in all
variants of the BC-method.
\smallskip

\noindent $\bullet$\quad We are grateful to the Referees for
useful remarks and criticism.

\section{Forward problem}\label{sec Forward problem}
To solve an inverse problem one has to study basic properties of
the solution to the corresponding forward problem. In this
section, we do it and derive the Duhamel-type representations for
the solution.

\subsection{Setup}\label{sec Setup}
Consider the initial boundary value problem of the form
 \begin{align}
\label{Dir 1} & iu_t+Ju_x+Vu=0, && x>0, \quad 0<t<T \\
\label{Dir 2} & u\big|_{t=0}=0, && x \geqslant 0 \\
\label{Dir 3} & u_1\big|_{x=0}=f, && 0\leqslant t \leqslant T\,,
 \end{align}
where $i=\sqrt{-1}$\,;
 \begin{equation*}
J:=\begin{pmatrix} 0&1\\
                  -1&0
   \end{pmatrix}, \qquad V=\begin{pmatrix} p&q\\
                  q&-p
   \end{pmatrix}\,,
 \end{equation*}
$V$ is a matrix {\it potential}, $p=p(x)$ and $q=q(x)$ being
real-valued functions (so that $V$ is symmetric and ${\rm
tr\,}V\equiv 0$); $f=f(t)$ is a complex-valued function ({\it
boundary control});
$u=u^f(x,t)=\begin{pmatrix} u_1(x,t)\\
u_2(x,t) \end{pmatrix}$ is a solution ({\it state}).

Define the domains on the plane
 $$ \Pi^T:=\{(x,t)\,|\,\,x\geqslant
0,\, 0\leqslant t \leqslant T\}, \,\,\, \Delta^\tau:=\{(x,t)\in
\Pi^T\,|\,\,x+T-\tau \leqslant t \}\,,
 $$
where $0<\tau \leqslant T$ is a parameter.
 \begin{theorem}\label{1}
Let $p,q \in C^1_{\rm loc}([0,\infty); \mathbb R)$. For any $f \in
C^1\left([0,T];{\mathbb C} \right)$ provided $f(0)=f^\prime(0)=0$,
problem (\ref{Dir 1})--(\ref{Dir 3}) has a unique classical
solution $u^f \in C^1\left(\Pi^T;{\mathbb C}^2 \right)$ such that
$u^f\big|_{t<x}=0$, and the representation
 \begin{equation}\label{repres u^f}
u^f(x,t)=u^f_*(x,t)+\int_x^t f(t-s)w(x,s)\,ds\,, \qquad x
\geqslant 0,\,\,0\leqslant t \leqslant T
 \end{equation}
holds, where $u^f_*(x,t):=f(t-x)\begin{pmatrix} 1\\
i \end{pmatrix}$ is the solution to (\ref{Dir 1})--(\ref{Dir 3}) for $V\equiv 0$, $w=\begin{pmatrix} w_1\\
w_2\end{pmatrix}$ is a vector-kernel such that $w\big|_{t<x}= 0$,
$w\big|_{\Delta^T} \in C^1(\Delta^T; {\mathbb C}^2)$, and
$w_1(0,\cdot)=0$.
 \end{theorem}
The proof will follow from the results of the next two sections.

\subsection{Solvability}\label{sec Solvability}
\subsubsection*{Auxiliary problem}
Consider the problem
 \begin{align}
\label{Aux 1} & iu_t+Ju_x=h, && x>0, \quad 0<t<T \\
\label{Aux 2} & u\big|_{t=0}=0, && x \geqslant 0 \\
\label{Aux 3} & u_1\big|_{x=0}=0, && 0\leqslant t \leqslant T
 \end{align}
with $h=\begin{pmatrix} h_1\\
h_2\end{pmatrix} \in C^1_{\rm loc}(\Pi^T;{\mathbb C}^2)$. Its
classical solution is obtained by integration as follows.
\smallskip

{\bf $1.$}\,\,\,Let $t>x$ and
$A(x,t),\,B(0,t-x),\,C(t-x,0),\,D(t+x,0)\in \Pi^T$ be the points
with the corresponding coordinates at the plane $\{(\xi,
\eta)\,|\,\,\xi {\geqslant} 0,\,\,\eta {\geqslant} 0\}$, $AB,\,BC,\,\dots$ the
straight lines connecting the points.
 \begin{align}
\label{u^1 for Aux t>x}
u_1(x,t)=-\frac{i}{2\sqrt{2}}\left[\int_{AB}-\int_{BC}+\int_{AD}\right]h_1\,
dl+\frac{1}{2\sqrt{2}}
\left[-\int_{AB}+\int_{BC}+\int_{AD}\right]h_2\,dl\,,\\
\label{u^2 for Aux t>x}
u_2(x,t)=-\frac{i}{2\sqrt{2}}\left[\int_{AB}-\int_{BC}+\int_{AD}\right]h_2\,
dl-\frac{1}{2\sqrt{2}}
\left[-\int_{AB}+\int_{BC}+\int_{AD}\right]h_1\,dl\,,
 \end{align}
where $dl$ is the length element on a curve in ${\mathbb R}^2$.
\smallskip

{\bf $2.$}\,\,\,Let $t\leqslant x$, and
$A(x,t),\,C(x-t,0),\,D(t+x,0)\in \Pi^T$. Then the representation
takes the form
 \begin{align}
\label{u^1 for Aux t<x}
u_1(x,t)=-\frac{i}{2\sqrt{2}}\left[\int_{AC}+\int_{AD}\right]h_1\,
dl+\frac{1}{2\sqrt{2}}
\left[-\int_{AC}+\int_{AD}\right]h_2\,dl\,,\\
\label{u^2 for Aux t<x}
u_2(x,t)=-\frac{i}{2\sqrt{2}}\left[\int_{AC}+\int_{AD}\right]h_2\,
dl-\frac{1}{2\sqrt{2}} \left[-\int_{AC}+\int_{AD}\right]h_1\,dl\,,
 \end{align}
As example, one has
 $$
\int_{AD}h_k\,dl=\sqrt{2}\int_0^t
h_k(-\eta+t+x,\eta)\,d\eta=\sqrt{2}\int_x^{t+x}
h_k(\xi,\xi+t-x)\,d\eta\,.
 $$
To justify (\ref{u^1 for Aux t>x})--(\ref{u^2 for Aux t<x}) is an
exercise in differentiation. The result can be written in the form
$u=Sh$ with an operator $S$ acting in $C^1_{\rm loc}(\Pi^T;
{\mathbb C}^2)$. It is easily seen that
 \begin{equation}\label{S}
{\rm supp\,}h \subset \Delta^\tau \quad \text{implies} \quad {\rm
supp\,}Sh \subset \Delta^\tau
 \end{equation}
for $0<\tau \leqslant T$, and
 \begin{equation}\label{Sh1=0 for x=0}
\left(Sh\right)_1\big|_{x=0}\,=\,0\,, \qquad 0 \leqslant t
\leqslant T\,.
 \end{equation}
The latter equality follows from behavior of the integration
lines: if $x \to 0$ then $A(x,t)$ and $B(0,t-x)$ tend to $A(0,t)$,
the line $AB$ vanishes, whereas the lines $AD$ and $BC$ coincide
in the limit. Therefore the sums  of the integrals in the square
brackets in (\ref{u^1 for Aux t>x}) tend to zero.

\subsubsection*{Integral equation}
Returning to problem (\ref{Dir 1})--(\ref{Dir 3}), we look for its
solution in the form
 \begin{equation}\label{u=dot u + w}
u^f=
u^f_*+w^f=f(t-x)\begin{pmatrix}1\\i\end{pmatrix}+\begin{pmatrix}w_1(x,t)\\w_2(x,t)\end{pmatrix}
 \end{equation}
and easily conclude that $w^f$ satisfies
\begin{align}
\label{w 1} & iw_t+Jw_x=-Vw-V u^f_*, && x>0, \quad 0<t<T \\
\label{w 2} & w\big|_{t=0}=0, && x \geqslant 0 \\
\label{w 3} & w_1\big|_{x=0}=0, && 0\leqslant t \leqslant T\,.
 \end{align}
Applying the operator $S$ and denoting $A:=-SV$, we arrive at the
equation
 \begin{equation}\label{w-Aw=A dot u}
w - Aw\,=\,A u^f_*
 \end{equation}
in $C^1_{\rm loc}(\Pi^T; {\mathbb C}^2)$, which is equivalent to
the original problem (\ref{Dir 1})--(\ref{Dir 3}). Its solvability
is established in the standard way. Namely,
\smallskip

\noindent $\bullet$\,\,\,since ${{\rm supp\,}} u^f_* \subset
\Delta^T$ holds, one has ${{\rm supp\,}}A u^f_*
\overset{(\ref{S})}\subset \Delta^T$ that allows to reduce
(\ref{w-Aw=A dot u}) to an equation in $C^1(\Delta^T; {\mathbb
C}^2)$
\smallskip

\noindent $\bullet$\,\,\,the space $C^1\left(\Delta^T ; {\mathbb
C}^2\right)$ is equipped with the family of semi-norms
 $$
\|y\|_\tau=\max\{\sup \limits_{\Delta^\tau}|y|,\,\sup
\limits_{\Delta^\tau}|\nabla y|\}\,, \qquad 0<\tau {\leqslant} T\,,
 $$
where $|y|:=\left[\sum_{k=1,2}|y_k|^2\right]^{\frac{1}{2}}$ and
$|\nabla y|:=\left[\sum_{k=1,2}|\nabla \Re y_k|^2+|\nabla \Im
y_k|^2\right]^{\frac{1}{2}}$, $\|y\|_T$ being a norm. As one can
easily derive from (\ref{u^1 for Aux t>x})--(\ref{u^2 for Aux
t<x}), the inequality
 $$
\|Ay\|_T\,\leqslant\,c \int_0^T \|y\|_\tau\,d\tau
 $$
is valid with a constant $c>0$. As a consequence, the operator
${\mathbb I}-A$ is an isomorphism in $C^1(\Delta^T ; {\mathbb
C}^2)$ and its inverse is represented via the Neumann series
$\left[{\mathbb I}-A\right]^{-1}=\sum_{j\geqslant 0}A^j$
converging in the operator norm: see, e.g., \cite{Nizh}, Theorem
4.1.1, page 198.
\smallskip

\noindent $\bullet$\,\,\,solving (\ref{w-Aw=A dot u}) by
$w=\sum_{j\geqslant 0}A^{j+1} u^f_*$ in $C^1(\Delta^T ; {\mathbb
C}^2)$ and extending the solution from $\Delta^T$ to $\Pi^T$ by
zero, we get the solution $w^f$ to (\ref{w 1})--(\ref{w 3}).
\smallskip

As a result, the right hand side of (\ref{u=dot u + w}) provides
the (unique) classical solution $u^f$ of the original problem
(\ref{Dir 1})--(\ref{Dir 3}) that proves the first part of Theorem
\ref{1}.

\subsection{Representation of solution}\label{sec Representation of solution}
\subsubsection*{Fundamental solution}
Fix $\varepsilon >0$; the function
 $$
\delta_\varepsilon(t):= \begin{cases} \frac{1}{\varepsilon}, &0\leqslant t \leqslant \varepsilon\\
0, & t>\varepsilon \end{cases}
 $$
tends to the Dirac $\delta$-function as $\varepsilon\to 0$. Put
$f=\delta_\varepsilon$ in (\ref{Dir 3}) and {\it define}
 \begin{equation}\label{u^delta eps}
u_*^{\delta_\varepsilon}(x,t):=
\delta_\varepsilon(t-x)\begin{pmatrix}1\\i\end{pmatrix}\,, \quad
u^{\delta_\varepsilon}:=
u_*^{\delta_\varepsilon}+w^{\delta_\varepsilon}\,,
 \end{equation}
where $w=w^{\delta_\varepsilon}$ satisfies
 \begin{equation}\label{w-Aw=A dot u delta eps}
w-Aw\,=\,A u_*^{\delta_\varepsilon} \qquad {\rm in}\,\,\,\Pi^T\,.
 \end{equation}
These definitions are quite meaningful. Indeed,
$u_*^{\delta_\varepsilon}$ is a piece-wise continuous ${\mathbb
C}^2$-valued function in $\Pi^T$ supported in the narrow strip
$\{(x,t)\in \Pi^T\,|\,\,x\geqslant 0,\,\,x \leqslant t \leqslant
x+\varepsilon\}$. As is seen from (\ref{u^1 for Aux
t>x})--(\ref{u^2 for Aux t<x}), operator $S$ (and, hence, operator
$A)$ can be applied to $u_*^{\delta_\varepsilon}$, whereas the
vector-function
 $$
A u_*^{\delta_\varepsilon}=-SV u_*^{\delta_\varepsilon}=-S\left[
(p+iq)\delta_\varepsilon
\begin{pmatrix} 1\\-i\end{pmatrix}\right]
 $$
is piece-wise continuous in $\Pi^T$ and supported in  $\Delta^T$.

If $\varepsilon \to 0$ then $u^{\delta_\varepsilon} \to
\delta_\varepsilon
\begin{pmatrix} 1\\i\end{pmatrix}=:u^\delta$ in the relevant sense of
distributions. For a fixed $(x,t)\in \Delta^T$, the value of $(A
u_*^{\delta_\varepsilon})(x,t)$ is expressed via the integrals in
(\ref{u^1 for Aux t>x}), which are taken over the small parts of
the contour $ABCD$ lying in the strip. Simple analysis shows that
$A u_*^{\delta_\varepsilon}$ converges in $L_\infty^{\rm
loc}(\Pi^T; {\mathbb C}^2)$ to a ${\mathbb C}^2$-valued function,
which we denote by $A u_*^{\delta}$. It is $C^1$-smooth in
$\Delta^T$ and vanishes for $t<x$.

The limit passage in (\ref{w-Aw=A dot u delta eps}) provides the
equation
 \begin{equation}\label{int eqn for w}
w-Aw=A u_*^\delta\,,
 \end{equation}
which is solvable by the arguments quite analogous to the ones
used above. Its solution $w=w^\delta(x,t)$ is $C^1$-smooth in
$\Delta^T$ and vanishes for $t<x$. Also, representing $w=-S[Vw+V
u_*^\delta]$, we arrive at the equality
$(w)_1\big|_{x=0}\overset{(\ref{Sh1=0 for x=0})}=0$, which is
declared in Theorem \ref{1}.

Thus, the limit passage in (\ref{u^delta eps}) leads to a
distribution of the form
 \begin{equation}\label{u^delta}
u^\delta(x,t):=[u_*^\delta+w^\delta](x,t)=\delta(t-x)\begin{pmatrix}
1\\i\end{pmatrix}+\begin{pmatrix} w_1(x,t)\\w_2(x,t)\end{pmatrix}
 \end{equation}
with ${\rm supp\,}w^\delta \subset \Delta^T$ and
$w^\delta\big|_{\Delta^T}\in C^1(\Delta^T; {\mathbb C}^2)$. It
corresponds to the control $f=\delta(t)$ in (\ref{Dir 3}) and is
called a {\it fundamental solution} to problem (\ref{Dir
1})--(\ref{Dir 3}).
\smallskip

The following result specifies the structure of the fundamental
solution.
\begin{lemma}\label{lemma p+iq=w2-iw1}
The equality
 \begin{equation}\label{p+iq=w2-iw1}
-iw_1(x,x+0)+w_2(x,x+0)=p(x)+iq(x)\,, \qquad x \geqslant 0
 \end{equation}
is valid.
\end{lemma}
$\square$\quad We use a relevant version of the WKB-method. Let
$k:=\begin{pmatrix} 1\\i\end{pmatrix}$ and
$\theta(s):=\frac{1}{2}\left[1+{\rm sign\,}s\right]$ be the
Heavyside function. Also, simplifying the notation, we write
$w^\delta\equiv w$ and $w(x,x+0)\equiv w(x,x)$. By
(\ref{u^delta}), one has
 $$
u^\delta(x,t)=\delta(t-x) k+\theta(t-x)w(x,t)\,, \qquad (x,t) \in
\Pi^T\,.
 $$
Differentiating $u^\delta$ as a distribution, we have
 \begin{align*}
& u^\delta_t=\delta'(t-x)
k+\delta(t-x)w(x,x)+\theta(t-x)w_t(x,t),\\
& u^\delta_x=-\delta'(t-x)k+
-\delta(t-x)w(x,x)+\theta(t-x)w_x(x,t),\\
& Vu^\delta=\delta(t-x)V(x)k+\theta(t-x)V(x)w(x,t)
 \end{align*}
that leads to
 \begin{align*}
& iu^\delta_t+Ju^\delta_x+Vu^\delta=\delta(t-x)\left[iw(x,x)-Jw(x,x)+V(x)k\right]+\\
& +\theta(t-x)\left[iw_t(x,t)+Jw_x(x,t)+V(x)w(x,t)\right]\,.
 \end{align*}
Since $u^\delta$ satisfies (\ref{Dir 1}) as a distribution, the
right hand side is equal to zero, and the expressions in the
square brackets have to vanish. In particular, one has
$iw(x,x)-Jw(x,x)+V(x)k=0$. Written by components, the latter
equality is equivalent to (\ref{p+iq=w2-iw1}). \quad
$\blacksquare$

\subsubsection*{Representation and generalized solutions}
As one can check, the classical solutions of (\ref{Dir
1})--(\ref{Dir 3}) are connected with $u^\delta$ via the Duhamel
formula
 \begin{equation}\label{Duhamel}
u^f(x,t)=\left[u^\delta(x, \cdot) \ast
f\right](t)\overset{(\ref{u^delta})}=f(t-x)\begin{pmatrix}
1\\i\end{pmatrix}+\int_x^t f(t-s)\begin{pmatrix}
w_1(x,s)\\w_2(x,s)\end{pmatrix}ds
 \end{equation}
that establishes (\ref{repres u^f}) and {\it completes the proof
of Theorem \ref{1}}.
\smallskip

Denote
 $$
\Omega^\xi\,:=\,\{x \in {\mathbb R}\,|\,\,0 \leqslant x \leqslant
\xi\}\,, \qquad \xi > 0\,.
 $$
In what follows we deal with controls $f\in L_2([0,T]; \mathbb
C)$. In this case, the right hand side in (\ref{Duhamel}) is well
defined as element of the space $C([0,T]; L_2(\Omega^T;{\mathbb
C}))$, and {\it is regarded} as a (generalized) solution of
(\ref{Dir 1})--(\ref{Dir 3}) for this class of controls.

The following general properties of solutions are easily derived
from the definitions and representation (\ref{Duhamel}): for $f\in
L_2\left([0,T]; \mathbb C\right)$ the relations
 \begin{equation}\label{supp u and delay relation}
{\rm supp\,}u^f(\cdot,\xi) \subset \Omega^\xi\,,\quad
u^{f_\xi}(\cdot,T)=u^f(\cdot,\xi)\,, \qquad 0< \xi \leqslant T
 \end{equation}
are valid, where $f_\xi(t):=f(t-(T-\xi))$ (recall Convention
\ref{Conv 1}!). The first relation shows that signals propagate
with velocity $1$, so that the finiteness of the domain of
influence does hold for the Dirac system. The second one follows
from the fact that the operator $J\frac{d}{dx} + V$, which governs
evolution of the system, does not depend on time.

One more specific property is the following. As was mentioned
above, equation (\ref{int eqn for w}), which determines the
regular part of the fundamental solution, is reduced to an
equation in $\Delta^T$. Therefore, $w$ is determined by the values
of the potential $V$  in $\Omega^T$ only (does not depend on
$V\big|_{x>T}$). One can easily see from (\ref{Duhamel}) that the
same is valid for the solution $u^f$: {\it it is determined by
$V\big|_{\Omega^T}$}.

\subsubsection*{Optimal setup}
Owing to the property mentioned above, the problem of the form
\begin{align}
\label{Dir 1 opt} & iu_t+Ju_x+Vu=0, &&  0<x<T,\,\, 0<t<T\\
\label{Dir 2 opt} & u\big|_{t<x}=0 \\
\label{Dir 3 opt} & u_1\big|_{x=0}=f, && 0\leqslant t \leqslant T
 \end{align}
turns out to be equivalent to (\ref{Dir 1})--(\ref{Dir 3}): to get
the solution of the latter problem one can solve (\ref{Dir 1
opt})--(\ref{Dir 3 opt}) and just extend it from
$\{(x,t)\,|\,\,x\in \Omega^T,\,\, t \in [0,T]\}$ to $\Pi^T$ by
zero. However, the form (\ref{Dir 1})--(\ref{Dir 3}) is optimal in
the sense that it does not invoke the values $V\big|_{x>T}$ which
have no influence on the solution.

\subsubsection*{Extended problem}
Denote $\Sigma^{2T}:=\{(x,t)\,|\,\,0<x<T, \, 0<t<2T-x\}$. Owing to
the FDI-property, problem (\ref{Dir 1 opt})--(\ref{Dir 3 opt})
possesses a natural extension of the form
\begin{align}
\label{Dir 1 ext} & iu_t+Ju_x+Vu=0, &&  (x,t)\in \Sigma^{2T}\\
\label{Dir 2 ext} & u\big|_{t<x}=0 \\
\label{Dir 3 ext} & u_1\big|_{x=0}=f, && 0\leqslant t \leqslant 2T
 \end{align}
that turns out to be well posed. Indeed, it can be reduced to a
relevant analog of equation (\ref{w-Aw=A dot u}) by the use of the
same formulas (\ref{u^1 for Aux t>x})--(\ref{u^2 for Aux t<x}), in
which the point $A(x,t)$ is permitted to take any position in
$\Sigma^{2T}$, whereas the integration is implemented over the
parts of the lines lying in $\Sigma^{2T}$. Also, repeating
`mutatis mutandis' the derivation of (\ref{Duhamel}), one gets its
analog of the Duhamel form
 \begin{equation}\label{Duhamel ext}
u^f(x,t)=f(t-x)\begin{pmatrix} 1\\i\end{pmatrix}+\int_x^t
f(t-s)\begin{pmatrix} w_1(x,s)\\w_2(x,s)\end{pmatrix}ds, \quad
(x,t)\in \Sigma^{2T}.
 \end{equation}
Problem (\ref{Dir 1 ext})--(\ref{Dir 3 ext}) is set up in the
domain $\Sigma^T$ and, hence, its solution $u^f$ {\it is
determined by} $V\big|_{\Omega^T}$. Along with the solution, a map
$f \mapsto u^f(0, \cdot)\big|_{0{\leqslant} t {\leqslant} 2T}$ is
also determined by the values of the potential in $\Omega^T$. The
latter fact will be used later.

The reason to regard problem (\ref{Dir 1 ext})--(\ref{Dir 3 ext})
as an extension of (\ref{Dir 1})--(\ref{Dir 3}) is the following.
Take a control $f$ in (\ref{Dir 3 ext}); let $\tilde
f:=f\big|_{[0,T]}$ be its restriction on the interval $[0,T]$.
Then, (\ref{Duhamel ext}) easily shows that for $t \leqslant T$,
the solution $u^f(\cdot,t)$ to (\ref{Dir 1 ext})--(\ref{Dir 3
ext}) coincides with the solution $u^{\tilde f}(\cdot,t)$ to
(\ref{Dir 1})--(\ref{Dir 3}):
\begin{equation}\label{u tilde f = u f}
u^{\tilde f}(\cdot,t)\,=\,u^f(\cdot,t) \qquad {\rm
in}\,\,\,\Omega^T\,,\,\,\,\,0 \leqslant t \leqslant T.
 \end{equation}

\section{Dynamical system}\label{sec Dynamical system}
\subsection{System $\alpha^T$}\label{sec System alpha T}
Here we consider problem (\ref{Dir 1 opt})--(\ref{Dir 3 opt}) as a
dynamical system, name it by $\alpha^T$, and endow with standard
attributes of control theory: spaces and operators.

\subsubsection*{Outer space}
The space ${{\mathcal F}^T}:=L_2\left([0,T]; {\mathbb C}\right)$
of controls with the inner product
 $$
(f,g)_{{{\mathcal F}^T}}:=\int_0^Tf(t)\bar g(t)\,dt
 $$
is said to be an {\it outer space} of system ${{{\alpha}^T}}$. It contains an
increasing family of subspaces
 $$
{\mathcal F}^{T,\xi}:=\{f \in {{\mathcal F}^T}\,|\,\,{\rm supp\,}f
\subset [T-\xi,T]\}, \qquad 0<\xi\leqslant T\,,
 $$
which consist of the delayed controls, $T-\xi$ and $\xi$ being a
delay and action time respectively. Also, there is a lineal set
 $$
{{\mathcal M}}^T:=\{f \in C^1\left([0,T]; {\mathbb
C}\right)\,|\,\,f(0)=f'(0)=0\}
 $$
of smooth controls producing classical solutions $u^f$. This set
is dense in ${{\mathcal F}^T}$.

\subsubsection*{Inner space}
The vector-function $u^f(\cdot,t)$ is a {\it state} of the system
at the moment $t$. The space of states  ${{\mathcal
H}^T}:=L_2(\Omega^T; {\mathbb C}^2)$ with the product
 $$
(a,b)_{{{\mathcal H}^T}} := \int_{\Omega^T}a(x) \cdot \bar
b(x)\,dx
 $$
($a\cdot \bar b:=a_1\bar b_1+ a_2\bar b_2$) is an {\it inner
space}. It contains an increasing family of subspaces
 $$
{{\mathcal H}^\xi}:=\{y \in {{\mathcal H}^T}\,|\,\,{\rm supp\,}y
\subset \Omega^\xi\}, \qquad 0<\xi \leqslant T\,.
 $$
By (\ref{supp u and delay relation}), one has
$u^f(\cdot,\xi)\subset {{\mathcal H}^\xi}$.

\subsubsection*{Control operator}
The input-state correspondence is realized by a {\it control
operator} $W^T_\alpha: {{\mathcal F}^T} \to {{\mathcal H}^T}$
 $${{W^T_{\alpha}}}  f\,:=\,
u^f(\cdot,T)\,.
 $$
This definition implies the operator relation on ${{\mathcal
M}}^T$ of the form
\begin{equation}\label{Wd/dt=LW alpha}
W^T_\alpha\left[-i\frac{d}{dt}\right]\,=\,L W^T_\alpha  \quad
\text{with} \quad L:=J\frac{d}{dx} +V\,,
 \end{equation}
which is just a way of writing the equation (\ref{Dir 1 opt}).
\smallskip

By (\ref{Duhamel}), one represents
 \begin{equation}\label{WAT repres}
(W^T_\alpha f)(x)=f(T-x)\begin{pmatrix} 1\\i\end{pmatrix}+\int_x^T
f(T-s)\begin{pmatrix} w_1(x,s)\\w_2(x,s)\end{pmatrix}ds
 \end{equation}
and see that ${{W^T_{\alpha}}} $ is a Fredholm operator. Its
integral part is a Volterra type operator; hence, $W^T_\alpha$
maps ${{\mathcal F}^T}$ {\it on its image} isomorphically.  By
(\ref{supp u and delay relation}), one has
 \begin{equation}\label{embedding in alpha^T}
W^T_\alpha{\mathcal F}^{T,\xi}\,\subset\, {{\mathcal
H}^\xi},\qquad 0<\xi\leqslant T\,.
 \end{equation}

\subsubsection*{Response operators}
$\bullet$\,\,The input-output correspondence is described by a
{\it response operator} $R^T_\alpha$, which acts in ${{\mathcal
F}^T}$ by the rule $R^T_\alpha: u^f_1(0,\cdot)\mapsto
u^f_2(0,\cdot)$, i.e.,
 $$
(R^T_\alpha f)(t)\,:=\,u^f_2(0,t), \quad 0 \leqslant t \leqslant
T\,.
 $$
Representation (\ref{Duhamel}) implies
 \begin{equation}\label{R^T repres}
(R^T_\alpha f)(t)=if(t)+\int_0^t r_\alpha(t-s) f(s)\,ds\,,\quad 0
\leqslant t \leqslant T\,,
 \end{equation}
where $r_\alpha:=w_2(0,\cdot) \in C^1([0,T]; \mathbb C)$ is a {\it
response function}. Hence, the response operator is bounded.
\smallskip

\noindent$\bullet$\, With system ${{{\alpha}^T}}$ one associates one more map
introduced via the extended problem (\ref{Dir 1 ext})--(\ref{Dir 3
ext}). It acts in the space $L_2([0,2T];{\mathbb C})$, is called
an {\it extended response operator} $R^{2T}_\alpha$ and defined by
 $$
(R^{2T}_\alpha f)(t)\,:=\,u^f_2(0,t), \quad 0 \leqslant t
\leqslant 2T\,.
 $$
Representation (\ref{Duhamel ext}) implies
\begin{equation}\label{R^2T repres}
(R^{2T}_\alpha f)(t)=if(t)+\int_0^t r_\alpha(t-s) f(s)\,ds\,,\quad
0 \leqslant t \leqslant 2T\,,
 \end{equation}
$r_\alpha:=w_2(0,\cdot) \in C^1([0,2T]; \mathbb C)$ is also called
a {\it response function}.

As is easy to recognize, owing to FDI, the operator
$R^{2T}_\alpha$ coincides with the (non-extended) response
operator of the system $\alpha^{2T}$. However, as was mentioned at
the end of sec \ref{sec Representation of solution}, operator
$R^{2T}_\alpha:f \mapsto u^f(0, \cdot)\big|_{0{\leqslant} t
{\leqslant} 2T}$ is determined by the potential
$V\big|_{\Omega^T}$ (does not depend on $V\big|_{x> T}$). It is
the fact, which motivates to regard $R^{2T}_\alpha$ as an
intrinsic object of system ${{{\alpha}^T}}$.

Note a simple relation
\begin{equation}\label{R^T via R^2T}
(R^T_\alpha f)(t)\,=\,\left(R^{2T}_\alpha f_T\right)(t+T)\,,
\qquad 0 \leqslant t \leqslant T
 \end{equation}
(here $f_T(t):=f(t-T)$), which easily follows from the delay
relation in (\ref{supp u and delay relation}). It shows that
$R^{2T}_\alpha$ determines $R^{T}_\alpha$.

To formulate one more property of $R^{2T}_\alpha$ we recall that,
for a linear operator $A$ acting from a space of ${\mathbb
C}^k$-valued functions to a space of ${\mathbb C}^l$-valued
functions, the map
\begin{equation}\label{def bar A}
\bar A: y\mapsto \overline{A \bar y}
 \end{equation}
is also a well-defined linear operator. Using the fact that
$R^{2T}_\alpha$ acts as a convolution with respect to time, one
can easily check the property
\begin{equation}\label{(RY)*=RY}
(Y^{2T}R^{2T}_\alpha)^*\,=\,Y^{2T}\bar R^{2T}_\alpha\,,
 \end{equation}
where $Y^{2T}$ acts in $L_2([0,2T];{\mathbb C})$ by the rule
 $$
(Y^{2T}f)(t):=f(2T-t), \qquad 0 \leqslant t \leqslant 2T\,.
 $$

\subsubsection*{Connecting operator}
A {\it connecting operator} $C^T_\alpha: {{\mathcal F}^T} \to
{{\mathcal F}^T}$ is
 $$
C^T_\alpha\,:=\,({{W^T_{\alpha}}} )^*\, {{W^T_{\alpha}}} \,.
 $$
Since ${{W^T_{\alpha}}} $ is an isomorphism on its image,
$C^T_\alpha$ turns out to be a bounded positive-definite operator
in ${{\mathcal F}^T}$. Its definition implies
\begin{equation}\label{(Cfg)=(uf ug) for alpha T}
\left(C^T_\alpha f,g\right)_{{{\mathcal F}^T}}=\left(W^T_\alpha f,
W^T_\alpha g\right)_{{{\mathcal H}^T}}=\left(u^f(\cdot,T),
u^g(\cdot,T)\right)_{{{\mathcal H}^T}}\,,
 \end{equation}
so that $C^T_\alpha$ connects the metrics of the outer and inner
spaces.

An important fact widely exploited in the BC-method is that the
connecting operator is determined by the response operator via
simple and explicit formula.
 \begin{lemma}\label{lemma CT via RT}
The representations
 \begin{align}
\notag & C^T_\alpha=-i\,[R^T_\alpha - (R^T_\alpha)^*], \qquad
(C^T_\alpha
f)(t)\,=\\
\label{CT via RT}& = 2f(t)-i \int_0^T\left[r_\alpha(t-s)-\bar
r_\alpha(s-t)\right]f(s)\,ds, \quad 0 \leqslant t \leqslant T
 \end{align}
are valid.
 \end{lemma}
$\square$\quad Take controls $f,g \in {{\mathcal M}}^T$ and define
a {\it Blagovestchenskii function}
 $$
b(s,t):=\left(u^f(\cdot,s), u^g(\cdot,t)\right)_{{{\mathcal
H}^T}}=\int_{\Omega^T}u^f(x,s) \cdot\, \bar u^g(x,t)\,dx, \,\,
(s,t)\in [0,T]\times[0,T]\,.
 $$
Since $u^f$ and $u^g$ are the classical solutions, (\ref{supp u
and delay relation}) implies $u^f\big|_{x=T}=u^g\big|_{x=T}=0$
that we use below. Also, we denote $L:=J
\partial_x+V$ and recall that $V$ is a real and symmetric
matrix, whereas $J$ is antisymmetric. Thereby, $L$ is symmetric by
Lagrange that we use in the equality $(\ast)$ below. Then, we have
  \begin{align}
\notag & i\,\left[b_s(s,t)+b_t(s,t)\right]=\int_{\Omega^T}\left[i
u^f_s(x,s)\cdot \bar u^g(x,t)-u^f(x,s)\overline{\left(i
u^g_t(x,t)\right)}\right]\,dx
\overset{(\ref{Dir 1})}=\\
\notag & =\int_{\Omega^T}\left[-\left(Lu^f\right) (x,s)\cdot \bar
u^g(x,t)+u^f(x,s)\cdot\overline{\left(L
u^g\right)(x,t)}\right]\,dx\overset{\ast}=\\
\notag & =\int_0^T\left[- Ju^f_x(x,s)\cdot \bar
u^g(x,t)+u^f(x,s)\cdot\overline{J u^g_x(x,t)}\right]\,dx=\\
\notag & =-\int_0^T\left[\left(J u^f(x,s)\right)_x \cdot\bar
u^g(x,t)+Ju^f(x,s)\cdot\bar
 u^g_x(x,t)\right]dx=\\
\notag & =-\left[Ju^f(x,s)\cdot\bar u^g(x,t)\right]\big|_{x=0}^{x=T}=Ju^f(0,s)\cdot \bar u^g(0,t)=\\
\notag & =u^f_2(0,s) \bar u^g(0,t)-u^f_1(0,s) \bar u^g_2(0,t)=
\left(R^T_\alpha f\right)(s) \bar g(t)- f(s)
\overline{\left(R^T_\alpha\right)
g(t)}=\\
\label{Calculation} & =:\, h(s,t)\,, \qquad (s,t)\in
[0,T]\times[0,T]\,.
  \end{align}
Integrating the equation $b_s+b_t=\frac{h}{i}$ for $s {\leqslant} t$ with
regard to $b(0,t)\overset{(\ref{Dir 2})}=0$, we easily get
 \begin{equation*}
b(s,t)=\frac{1}{i}\int_{t-s}^t h(\eta-t+s, \eta)\,d\eta
 \end{equation*}
that leads to
 \begin{align}
\notag & b(T,T)=\frac{1}{i}\int_0^T h(\eta,
\eta)\,d\eta=\frac{1}{i} \int_0^T\left[(R^T_\alpha f)(\eta) \bar
g(\eta)- f(\eta)
\overline{(R^T_\alpha g(\eta)}\right]\,d\eta=\\
\label{*} &=\left(-i\left[R^T_\alpha -
\left(R^T_\alpha\right)^*\right] f, g\right)_{{{\mathcal F}^T}}\,.
 \end{align}
On the other hand, the definition of $b$ implies
\begin{equation}\label{**}
b(T,T)=\left(u^f(\cdot,T), u^g(\cdot,T)\right)_{{{\mathcal
H}^T}}\overset{(\ref{(Cfg)=(uf ug) for alpha T})}=\left(C^T_\alpha
f,g\right)_{{{\mathcal F}^T}}\,.
 \end{equation}
Comparing (\ref{*}) with (\ref{**}), taking into account
arbitrariness of $f,g$ and density of ${{\mathcal M}}^T$ in
${{\mathcal F}^T}$, we arrive at the first of representations
(\ref{CT via RT}).

Substituting (\ref{R^T repres}) to the first representation, one
easily gets the second one, the kernel in the square brackets
being written in accordance with Convention \ref{Conv 1}. \quad
$\blacksquare$

\subsubsection*{Controllability}
The sets
 $$
{{\mathcal U}}^\xi_\alpha\,:=\,W^T_\alpha {\mathcal F}^{T,\xi}
\overset{(\ref{supp u and delay
relation})}=\{u^f(\cdot,\xi)\,|\,\,f \in {{\mathcal F}^T}\}\,,
\qquad 0<\xi \leqslant T
 $$
are called {\it reachable}. Since $W^T_\alpha$ is a Fredholm
operator (see (\ref{repres u^f})), each ${{\mathcal U}}^\xi_\alpha$ is a
closed subspace in ${{\mathcal H}^T}$. By (\ref{embedding in
alpha^T}), the embedding ${{\mathcal U}}^\xi_\alpha \subset {{\mathcal H}^\xi}$ occurs. In
the mean time, an important fact is that this embedding {\it is
not dense}.

Indeed, consider the case $\xi=T$ and choose a nonzero $h \in
{{\mathcal U}}^T_\alpha$. Then there is an $f \in {{\mathcal F}^T}$ such that
$u^f(\cdot,T)=h$. By (\ref{Duhamel}), the latter equality takes
the form
 $$
\begin{cases}
(A^T f)(x):=f(T-x)+ \int_x^Tw_1(x,s) f(T-s)\,ds = h_1(x)\\
(B^T f)(x):=if(T-x)+ \int_x^Tw_2(x,s) f(T-s)\,ds = h_2(x)
\end{cases},\qquad  x \in \Omega^T
 $$
of two 2-nd kind Volterra equations. Solving any of them, one
determines $f$ uniquely. Therefore, there is an isomorphism
$M^T=B^T(A^T)^{-1}$ in $L_2(\Omega^T; {\mathbb C})$, which links
the components: the relation $h_2=M^T h_1$ is necessary and
sufficient for $h$ to belong to ${{\mathcal U}}^T_\alpha$. As a
consequence, taking a nonzero $a \in {{\mathcal H}^T}$ provided
$a_1=-(M^T)^*a_2$, one gets
 \begin{align*}
&\left(u^f(\cdot,T),a\right)_{{{\mathcal
H}^T}}=\left(\begin{pmatrix}u^f_1(\cdot,T)\\M^T
u^f_1(\cdot,T)\end{pmatrix},\begin{pmatrix}-(M^T)^*a_2\\a_2\end{pmatrix}\right)_{{{\mathcal H}^T}}=\\
&=\left(u^f_1(\cdot,T),-(M^T)^*a_2\right)_{L_2(\Omega^T; {\mathbb
C})}+\left(M^T u^f_1(\cdot,T),a_2\right)_{L_2(\Omega^T; {\mathbb
C})}=0\,,
 \end{align*}
for any $f$, i.e., $a \in {{\mathcal H}^T}
\ominus{{\mathcal U}}^T_\alpha\not=\{0\}$. The same is valid for all $\xi>0$:
one has ${\mathcal H}^\xi \ominus{{\mathcal U}}^\xi_\alpha\not=\{0\}$.

Thus, at any moment $t=\xi$ the states $u^f(\cdot, \xi)$ fill up
the segment $\Omega^\xi$ (see (\ref{supp u and delay relation}));
they belong to the subspace ${\mathcal H}^\xi$ but do not constitute a
complete system in it. In terms of control theory, this fact is
interpreted as {\it a lack of local boundary controllability}. By
philosophy of the BC-method, such a lack ever implies
complications and difficulties in solving the corresponding
inverse problem. That is why, in the next section we construct an
auxiliary system $\beta^T$, which does not possess such a
`drawback'.

\subsection{System $\beta^T$}\label{sec System beta T}
\subsubsection*{System $\bar\alpha^T$}
This system is associated with the problem
\begin{align}
\label{Dir 1 opt v} & iv_t-Jv_x-Vv=0, &&  0<x<T,\,\, 0<t<T\\
\label{Dir 2 opt v} & v\big|_{t<x}=0 \\
\label{Dir 3 opt v} & v_1\big|_{x=0}=g, && 0\leqslant t \leqslant
T
 \end{align}
and its extension by FDI
\begin{align}
\label{Dir 1 ext v} & iv_t-Jv_x-Vv=0, &&  (x,t)\in \Sigma^{2T}\\
\label{Dir 2 ext v} & v\big|_{t<x}=0 \\
\label{Dir 3 ext v} & v_1\big|_{x=0}=g, && 0\leqslant t \leqslant
2T\,.
 \end{align}
Both of the problems are well posed; their solutions $v=v^g(x,t)$
are connected with the solutions to (\ref{Dir 1 opt})--(\ref{Dir 3
opt}) and (\ref{Dir 1 ext})--(\ref{Dir 3 ext}) in a simple way:
 \begin{equation}\label{v = bar bar u}
v^g(x,t)\,=\,\overline{u^{\bar g}(x,t)}\,.
 \end{equation}
Therefore, the attributes of ${{{\bar \alpha}^T}}$ are simply connected with the
ones of ${{{\alpha}^T}}$ as follows. Recall the notation (\ref{def bar A}).
\smallskip

\noindent$\bullet$\,\,The outer and inner spaces of ${{{\bar \alpha}^T}}$ are the
same, i.e., ${{\mathcal F}^T}$ and ${{\mathcal H}^T}$.
\smallskip

\noindent$\bullet$\,\,The control operator ${{W^T_{\bar\alpha}}}: {{\mathcal F}^T}
\to {{\mathcal H}^T},\,\,{{W^T_{\bar\alpha}}} g:=v^g(\cdot,T)$ satisfies
${{W^T_{\bar\alpha}}}\overset{(\ref{v = bar bar u})}=\bar W^T_\alpha$ an is
represented in the form
\begin{equation}\label{WATb repres}
(W^T_{\bar \alpha} g)(x)\overset{(\ref{WAT
repres})}=g(T-x)\begin{pmatrix} 1\\-i\end{pmatrix}+\int_x^T
g(T-s)\begin{pmatrix} \bar w_1(x,s)\\\bar
w_2(x,s)\end{pmatrix}ds\,.
 \end{equation}
By analogy with (\ref{Wd/dt=LW alpha}), the definition of
$W^T_{\bar \alpha}$ implies the relation on ${{\mathcal M}}^T$:
 \begin{equation}\label{Wd/dt=LW bar alpha}
W^T_{\bar \alpha}\left[-i\frac{d}{dt}\right]\,=\,-\,L W^T_{\bar
\alpha}\,.
 \end{equation}
\smallskip

\noindent$\bullet$\,\,The response operators
 $$
R^{T}_{\bar \alpha} g:=v^g_2(0,\cdot)\big|_{0{\leqslant} t
{\leqslant} T}={\bar R}^T_\alpha g \quad \text{and} \quad
R^{2T}_{\bar \alpha} g:=v^g_2(0,\cdot)\big|_{0{\leqslant} t
{\leqslant} 2T}={\bar R}^{2T}_\alpha g
 $$
act in ${{\mathcal F}^T}$ and ${{\mathcal F}^{2T}}$ respectively.
\smallskip

\noindent$\bullet$\,\,The connecting operator $C^T_{\bar \alpha}:
{{\mathcal F}^T} \to {{\mathcal F}^T}$ is defined by $C^T_{\bar
\alpha}:=({{W^T_{\bar\alpha}}})^*{{W^T_{\bar\alpha}}}$ and admits
the representations
\begin{align}
\notag & C^T_{\bar \alpha}\,=\,i\,[R^T_{\bar \alpha} - (R^T_{\bar
\alpha})^*]\,= \,i\,[\bar R^T_\alpha - (\bar R^T_\alpha)^*], \quad
(C^T_\alpha
f)(t)\,=\\
\label{CTb via RTb}& = 2f(t)+i \int_0^T\left[\bar r_\alpha(t-s)-
r_\alpha(s-t)\right]f(s)\,ds, \quad 0 \leqslant t \leqslant T \,,
 \end{align}
which are derived quite analogously to (\ref{CT via RT}).
\smallskip

In regard to controllability, it is easy to recognize that system
${{{\bar \alpha}^T}}$ does possess the same `drawback': for its reachable sets
${{\mathcal U}}^\xi_{\bar \alpha}:={{W^T_{\bar\alpha}}} {{\mathcal F}^{T,\xi}}$ one has ${\mathcal
H}^\xi \ominus{{\mathcal U}}^\xi_{\bar \alpha}\not=\{0\}$, i.e., the system
is not locally controllable from the endpoint $x=0$.

\subsubsection*{Restoring controllability}
Now, introduce a dynamical system ${{{\beta}^T}}$ with the following
attributes.
\smallskip

\noindent$\bullet$\,\,The outer space is ${\mathcal F}^T_\beta:={{\mathcal
F}^T} \oplus {{\mathcal F}^T}=L_2([0,T];{\mathbb C}^2)$; its elements
are represented as columns
 \begin{equation}\label{repres matrix controls}
  \begin{pmatrix}f\\g\end{pmatrix} \in {{{\mathcal F}^T_{\beta}}}\,, \qquad f,g \in {{\mathcal F}^T}\,.
 \end{equation}
The outer space contains the delayed subspaces ${\mathcal
F}^{T,\xi}_\beta:={{\mathcal F}^{T,\xi}} \oplus {{\mathcal F}^{T,\xi}}$
and the smooth class ${{\mathcal M}}^T_\beta:={{\mathcal M}}^T \oplus {{\mathcal M}}^T$.
\smallskip

\noindent$\bullet$\,\,The inner space is ${{\mathcal H}^T}$. Recall
that it contains the subspaces ${{\mathcal H}^\xi}=\{y \in {{\mathcal
H}^T}\,|\,\,{\rm supp\,}y \subset
\Omega^\xi\}\,\,\,\,(0<\xi\leqslant T)$.
\smallskip

\noindent$\bullet$\,\,In the matrix terms corresponding to
(\ref{repres matrix controls}), the control operator $W^T:{{\mathcal
F}^T}\to {{\mathcal H}^T}$ is defined as an operator row
$W^T:=(W^T_\alpha\,\,\,W^T_{\bar \alpha})$, so that
\begin{equation}\label{def W beta}
W^T\begin{pmatrix}f\\g\end{pmatrix}=W^T_\alpha f + W^T_{\bar
\alpha} g=u^f(\cdot,T)+v^g(\cdot,T)\,.
 \end{equation}
The definition implies the operator relation
\begin{equation}\label{Wd/dt=LW}
W^T\left[-i Q \frac{d}{dt}\right]\,=\,L W^T
 \end{equation}
with $L=J\frac{d}{dx} +V$ and
$Q:=\begin{pmatrix}1&0\\0&-1\end{pmatrix}$, which just summarizes
(\ref{Wd/dt=LW alpha}) and (\ref{Wd/dt=LW bar alpha}). Also, in
view of
 $$
W^T\begin{pmatrix}f\\0\end{pmatrix}=W^T_\alpha f = u^f(\cdot,
T)\,,\quad W^T\begin{pmatrix}0\\g\end{pmatrix}=W^T_{\bar \alpha} g
= v^g(\cdot, T)\,,
 $$
it is reasonable to regard ${{{\alpha}^T}}$ and $\bar \alpha^T$ as the
subsystems of $\beta^T$.

Denote
\begin{equation}\label{def varkappa and check w}
\varkappa :=\begin{pmatrix} 1&1\\i&-i\end{pmatrix}, \quad \check
w(x,s):= \begin{pmatrix} w_1(x,s)&\bar w_1(x,s)\\w_2(x,s)&\bar
w_2(x,s)\end{pmatrix},
 \end{equation}
and define $Y^T: {{{\mathcal F}^T_{\beta}}} \to {{\mathcal H}^T},\,\,\,(Y^T
a)(x):=a(T-x),\,\,\,x \in \Omega^T$.  By (\ref{WAT repres}) and
(\ref{WATb repres}), for an $a=\begin{pmatrix}f\\g\end{pmatrix}\in
{{{\mathcal F}^T_{\beta}}}$ one has
 \begin{align}
\notag & (W^T a)(x)=\varkappa a(T-x)+\int_x^T
\check w(x,s)a(T-s)\,ds=\\
\label{WT repres}&  =: \left(\varkappa\left[{\mathbb
I}+M^T\right]Y^T a\right)(x)\,, \qquad x \in \Omega^T\,.
 \end{align}
In this representation, $\varkappa Y^T$ is an isomorphism from
${\mathcal F}^T_\beta$ onto ${{\mathcal H}^T}$ and $M^T$ is a Volterra
operator in ${{\mathcal H}^T}$. Hence, $W^T$ maps ${\mathcal F}^T_\beta$
{\it onto} ${{\mathcal H}^T}$ isomorphically. By (\ref{supp u and
delay relation}), the same holds for intermediate $\xi$'s and we
arrive at the relations
 \begin{equation}\label{U^xi=H^xi in beta}
{{\mathcal U}}^\xi_\beta:= W^T {\mathcal F}^{T,\xi}_\beta\,=\,{{\mathcal
H}^\xi}\,, \qquad 0<\xi \leqslant T\,.
 \end{equation}

Thus, in contrast to $\alpha^T$ and ${\bar \alpha}^T$, system
$\beta^T$ does possess the property of local exact
controllability: its reachable sets exhaust the (sub)spaces
${{\mathcal H}^\xi}$, which they belong to.

\subsubsection*{Connecting operator}
In system $\beta^T$, its connecting operator acts in ${\mathcal
F}^T_\beta$ by the rule
 \begin{align}
\notag & C^T:=(W^T)^*W^T =\begin{pmatrix}
(W^T_{\alpha})^*\\(W^T_{\bar
\alpha})^*\end{pmatrix}\begin{pmatrix} W^T_{\alpha} \,\, W^T_{\bar
\alpha}\end{pmatrix}=\begin{pmatrix}
(W^T_{\alpha})^*W^T_{\alpha}&(W^T_{\alpha})^*W^T_{\bar
\alpha}\\(W^T_{\bar \alpha})^*W^T_{\alpha}&(W^T_{\bar
\alpha})^*W^T_{\bar
\alpha}\end{pmatrix}=\\
\notag & =\begin{pmatrix} C^T_{\alpha \alpha} & C^T_{\alpha \bar
\alpha}\\C^T_{\bar \alpha \alpha}&C^T_{\bar \alpha \bar
\alpha}\end{pmatrix}.
 \end{align}
Since $W^T$ is an isomorphism, $C^T$ turns out to be a {\it
positive definite isomorphism} in ${{\mathcal F}^T}_\beta$.

For our approach to inverse problems, the following fact is of
crucial character: operator $C^T$ is determined by the response
operator $R^{2T}_\alpha$ of the (sub)system $\alpha^T$. To
formulate the result let us introduce an extension by zero $Z^T:
{\mathcal F}^T_\beta\to L_2\left([0,2T];{\mathbb C}\right)$,
 $$
(Z^Tf)(t)\,:=\,\begin{cases} f(t), &0 \leqslant t \leqslant T\\0,
& T< t \leqslant 2T
                          \end{cases}\,.
 $$
Its conjugate $(Z^T)^*$ restricts controls $f \in L_2\left([0,2T];
{\mathbb C}\right)$ on the segment $[0,T]$. Also, recall
Convention \ref{Conv 1}.
\begin{lemma}\label{lemma CT beta via R2T alpha}
The representations
 \begin{align}
\notag & C^T\,=\,
\begin{pmatrix} -i\,[R^T_\alpha -
(R^T_\alpha)^*] &
-i(Z^T)^* Y^{2T}\bar R^{2T}_\alpha Z^T\\
i\,(Z^T)^* Y^{2T} R^{2T}_\alpha Z^T & i\,[\bar R^T_\alpha -
(\bar R^T_\alpha)^*]\end{pmatrix}\,,\\
\label{CT beta via R2T alpha} & \left(C^T a\right)(t)=2
a(t)+\int_0^Tc^T(t,s) a(s)\,ds\,, \qquad 0 \leqslant t \leqslant T
 \end{align}
are valid, where $c^T(t,s)$ is a matrix kernel with the elements
 \begin{align}
\notag & c_{11}(t,s)(t)=-i\,[r_\alpha(t-s)-\bar r_\alpha(s-t)]\,,
\quad
c_{12}(t,s)=-i\,\bar r_\alpha(2T-t-s)\,,\\
& c_{21}(t,s)= i\,r_\alpha(2T-t-s)\,,\quad c_{22}(t,s)(t)=i\,[\bar
r_\alpha(t-s)-r_\alpha(s-t)]\label{kernel c in detail in forward
problem}\,,
 \end{align}
which are $C^1$-smooth outside the diagonal $t=s$.
\end{lemma}
$\square$

\noindent $\bullet$\,\,For the diagonal blocks $C^T_{\alpha
\alpha}=C^T_\alpha$ and $C^T_{\bar \alpha \bar \alpha}=C^T_{\bar
\alpha}$, the representation does hold by (\ref{CT via RT}) and
(\ref{CTb via RTb}).
\smallskip

\noindent $\bullet$\,\,Choose $f \in C_0^\infty\left((0,T);
\mathbb C\right)$ and $g \in {{\mathcal M}}^T$. Denote $f_0:=Z^Tf
\in\\
C_0^\infty\left((0,2T); \mathbb C\right))$; let $u^{f_0}$ and
$v^g$ be the (classical) solutions to problems (\ref{Dir 1
ext})--(\ref{Dir 3 ext}) and (\ref{Dir 1 opt v})--(\ref{Dir 3 opt
v}) respectively. A Blagovestchenskii function of the form
 $$
b(s,t):=\left(u^{f_0}(\cdot, s), v^g(\cdot,t)\right)_{{{\mathcal
H}^T}}=\int_{\Omega^T}u^{f_0}(x, s)\, \overline{v^g(x,t)}\,dx
 $$
is evidently well defined for $s,t \leqslant T$.  Let $s$ and $t$
be such that $T\leqslant s \leqslant 2T$ and $t\leqslant 2T-s$. In
this case, $u^{f_0}(\cdot,s)$ is defined on the segment
$0\leqslant x \leqslant 2T-s$, which is bigger than the segment
$0\leqslant x \leqslant t$, where $v^g(\cdot,t)$ is supported: see
(\ref{v = bar bar u}) and (\ref{supp u and delay relation}).
Therefore, the product $u^{f_0}(x,s)v^g(x,t)$ is well defined and
supported on the smaller segment $0\leqslant x \leqslant t$, and
can be extended to $t< x \leqslant T$ by zero as a $C^1$-smooth
function of $x$. Under this agreement, $b$ is a well-defined
$C^1$-smooth function in the domain
$\Theta^T:=\{(s,t)\,|\,\,0\leqslant s \leqslant 2T,\,\,\,
t\leqslant \min\{T,2T-s\}\}$.

In perfect analogy to the calculations in (\ref{Calculation}), one
can derive an equation
 \begin{equation*}
i\,[b_t-b_s]\,=\,h \qquad \text{in}\,\,\Theta^T
 \end{equation*}
with  $h(s,t):=-\left(R^{2T}_\alpha f_0\right)(s)\bar
g(t)+f_0(s)\overline{\left(R^T_{\bar \alpha}g\right)(t)}$.
Integrating this equation with regard to
$b\big|_{t=0}\overset{(\ref{Dir 2 opt v})}=0$, one gets
 $$
b(s,t)\,=\,-i\int_0^t h(s+t-\eta, \eta)\,d \eta
 $$
that implies
 \begin{align}
\notag & b(T,T)=i \int_0^T\left[\left(R^{2T}_\alpha
f_0\right)(2T-\eta)\bar
g(\eta)-f_0(2T-\eta)\overline{\left(R^T_{\bar
\alpha}g\right)(\eta)}\right]\,d\eta=\\
\label{++}& =\int_0^T\left(R^{2T}_\alpha f_0\right)(2T-\eta)\bar
g(\eta)\,d\eta=\left(i(Z^T)^*Y^{2T}R^{2T}_\alpha Z^T f,
g\right)_{{{\mathcal F}^T}}
 \end{align}
since $f_0(2T-\eta)\big|_{0{\leqslant} \eta {\leqslant} T}\equiv
0$.
\smallskip

\noindent $\bullet$\,\,\,By the definition of $b$, one has
 \begin{align}
\notag & b(T,T)=\left(u^{f_0}(\cdot, T),
v^g(\cdot,T)\right)_{{{\mathcal H}^T}}\overset{(\ref{u tilde f = u
f})}=\left(u^f(\cdot, T), v^g(\cdot,T)\right)_{{{\mathcal H}^T}}=\\
\label{+++} & = \left(W^T_\alpha f, W^T_{\bar
\alpha}g\right)_{{{\mathcal H}^T}}=\left(C^T_{\bar \alpha \alpha} f,
g\right)_{{{\mathcal F}^T}}\,.
 \end{align}
Comparing (\ref{++}) with (\ref{+++}) and taking into account the
arbitrariness of $f,g$, one gets
 $$
C^T_{\bar \alpha \alpha}\,=\,i(Z^T)^*\,Y^{2T}R^{2T}_\alpha Z^T\,.
 $$
Since $C^T$ is self-adjoint in ${{\mathcal F}^T}_\beta$, one arrives
at
$$
C^T_{\alpha \bar \alpha}=\left(C^T_{\bar \alpha
\alpha}\right)^*\,=\,-i(Z^T)^*\left[Y^{2T}R^{2T}_\alpha\right]^*
Z^T\overset{(\ref{(RY)*=RY})}=-i(Z^T)^*\,Y^{2T}\bar R^{2T}_\alpha
Z^T
 $$
that proves the first representation in (\ref{CT beta via R2T
alpha}).
\smallskip

\noindent $\bullet$\,\,\,The expressions for $c_{kl}$ follow from
(\ref{R^T repres}) and (\ref{R^2T repres}). Since $r_\alpha \in
C^1\left([0,2T]; \mathbb C\right)$, the elements $c_{12}$ and
$c_{21}$ are $C^1$-smooth, whereas $c_{11}$ and $c_{22}$ are
$C^1$-smooth outside the diagonal. Note that the non-diagonal
terms $-i(Z^T)^* Y^{2T}\bar R^{2T}_\alpha Z^T$ and $i(Z^T)^*
Y^{2T}R^{2T}_\alpha Z^T$ do not contribute to the `leading part'
$2\mathbb I$ of $C^T$ owing to the identity
$(Z^T)^*Y^{2T}Z^T=\mathbb O$. \quad $\blacksquare$
\smallskip

Note an important fact: (\ref{CT beta via R2T alpha}) shows that
the connecting operator $C^T$ of the system $\beta^T$ {\it is
fully determined} by the response operator $R^{2T}_\alpha$ of the
basic system $\alpha^T$.

\section{Inverse problem}\label{sec Inverse problem}
\subsection{Main result}\label{subsec Main result}
As was more than once mentioned above, all the objects of system
${{{\alpha}^T}}$ are determined by the potential $V$ in $\Omega^T$
(do not depend on $V\big|_{x>T}$). In particular, this holds for
the response operator $R^{2T}_\alpha$. Such a character of
dependence motivates the following setup of the inverse problem
(IP):

\noindent {\it for an (arbitrary) fixed $T>0$, given $R^{2T}$ to
recover $V\big|_{\Omega^T}$.}

\noindent It is the so-called time-optimal setup, which is most
relevant to the FDI-property of the dynamical Dirac system.

By (\ref{R^2T repres}), to give $R^{2T}_\alpha$ is to know the
response function $r_\alpha\big|_{0 {\leqslant} t {\leqslant} 2T}$, which may be
regarded as the {\it inverse data}. To solve IP is
\smallskip

\noindent (i)\,\,\, to provide a procedure that recovers the
potential via the response function and
\smallskip

\noindent (ii)\,\,\,to characterize the inverse data that means to
describe the necessary and sufficient conditions on $r_\alpha$,
which provide solvability of the IP.
\smallskip

\noindent We do it in sec \ref{subsec Solving IP}. However, the
main result of the paper concerns to the item (ii) and is the
following. Let us specify again that we deal with the systems of
the form (\ref{Dir 1 opt})--(\ref{Dir 3 opt}) with a real
symmetric zero-trace potentials $V$.
\smallskip

Let $r=r(t)$ be a $\mathbb C$-valued function given on
$0{\leqslant} t {\leqslant} 2T$. With $r$ one associates

\noindent$\bullet$\,\,\,an operator $R$ acting in
$L_2\left([0,2T];\mathbb C\right)$ by the rule
 \begin{equation}\label{cal R repres}
\left(R f\right)(t):=if(t)+\int_0^t r(t-s) f(s)\,ds\,,\qquad 0
\leqslant t \leqslant 2T
 \end{equation}

\noindent$\bullet$\,\,\,an operator $R_T$ acting in
$L_2\left([0,T];\mathbb C\right)$ by the rule
 \begin{equation}\label{cal RT repres}
\left(R_T f\right)(t):=\left(R f_T\right)(t+T)\,,\qquad 0
\leqslant t \leqslant T\,,
 \end{equation}
where $f_T(t)=f(t-T)$

\noindent$\bullet$\,\,\,an operator $C$ acting in
$L_2([0,2T];{\mathbb C}^2)=L_2([0,2T];{\mathbb C})\oplus
L_2([0,2T];{\mathbb C})$ by the rule
 \begin{equation}\label{cal C repres} C\,:=\,
\begin{pmatrix} -i\,[R_T -
(R_T)^*] & -i\,(Z^T)^* Y^{2T}\bar R Z^T\\i\,(Z^T)^* Y^{2T} R Z^T &
i\,[\bar R_T - (\bar R_T)^*] \end{pmatrix}
 \end{equation}
(recall (\ref{def bar A})). As one can easily recognize, these
definitions are inspired by the rep\-re\-sen\-tations (\ref{R^2T
repres}), (\ref{R^T via R^2T}), and (\ref{CT beta via R2T alpha})
respectively.
 \begin{theorem}\label{Th Charact}
A function $r\in C^1\left([0,2T]; \mathbb C\right)$ is the
response function of a system  $\alpha^T$ with a $C^1$-smooth
potential $V$ if and only if the operator $C$ is a positive
definite isomorphism.
 \end{theorem}
The proof is postponed to sec \ref{subsec Solving IP}. In the next
section we provide the devices for the procedure that recovers the
potential.

\subsection{A portion of Operator Theory}\label{subsec A portion of Operator Theory}

The results of this section are of general character. An abstract
operator scheme, which is developed below, is the core of the
one-dimensional BC-method \cite{1D-BC, BV_JIPP_2}. We use the
notation, which prepares its application to the system $\beta^T$.

\subsubsection*{Lemmas on projectors}

\begin{lemma}\label{lemma cPxi via blocks}
Let

\noindent$\bullet$\,\, ${\mathcal F}$ and ${\mathcal F}'\subset {\mathcal F}$
be a Hilbert space and its subspace

\noindent$\bullet$\,\,$C$ a positive definite isomorphism in
${\mathcal F}$

\noindent$\bullet$\,\,${{\mathcal P}}$ the projector in ${\mathcal F}$ onto
${\mathcal F}'$ in parallel to the subspace $C^{-1}[{\mathcal F} \ominus
{\mathcal F}']$.

\noindent Then the representation and equality
 \begin{equation}\label{repres cP via blocks}
{{\mathcal P}}=e\,\left[e^* C e\right]^{-1}e^* C \quad \text{and} \quad
C{{\mathcal P}}\,=\,{{\mathcal P}}^*C
 \end{equation}
are valid, where $e: {\mathcal F}' \to {\mathcal F}$ is the (natural)
embedding.

Let ${\mathcal F}_C$ be ${\mathcal F}$ endowed with the inner product
$(a,b)_{{\mathcal F}_C}:=(Ca,b)_{\mathcal F}$. Then ${{\mathcal P}}$ is the
orthogonal projector in ${\mathcal F}_C$ onto ${\mathcal F}'$.
 \end{lemma}
$\square$\quad Since $C$ is an isomorphism, the projector ${{\mathcal
P}}$ is a well-defined bounded operator. The block $e^* C e$ of
$C$ in ${\mathcal F}'$ is a positive definite isomorphism, so that its
inverse $[e^* C e]^{-1}$ is well defined.

The right hand side of the representation in (\ref{repres cP via
blocks}) obeys the characteristic properties of the projector in
question: ${{\mathcal P}}^2={{\mathcal P}},\,\,\,{{\mathcal P}} {\mathcal F}={\mathcal
F}'$, and ${{\mathcal P}} C^{-1}[{\mathcal F} \ominus {\mathcal F}']=\{0\}$,
the latter following from the evident $e^* [{\mathcal F} \ominus {\mathcal
F}']=\{0\}$. Hence, $[e^* C e]^{-1}e^* C={{\mathcal P}}$. The equality
in (\ref{repres cP via blocks}) easily follows from the
representation.

By the equality, one has
 $$
({{\mathcal P}} a, b)_{{\mathcal F}_C}=(C{{\mathcal P}} a, b)_{{\mathcal F}}
=({{\mathcal P}}^*Ca, b)_{{\mathcal F}}=(Ca, {{\mathcal P}} b)_{{\mathcal
F}}=(a,{{\mathcal P}} b)_{{\mathcal F}_C},
 $$
so that ${{\mathcal P}}^2={{\mathcal P}},\,\,\,{{\mathcal P}} {\mathcal F}={\mathcal F}'$, and ${{\mathcal P}}$ is
self-adjoint in ${\mathcal F}_C$. Hence, it is the orthogonal
projector in ${\mathcal F}_C$ onto ${\mathcal F}'$. \quad $\blacksquare$
\smallskip

Now, let ${\mathcal F}=L_2\left([0,T]; {\mathbb C}^n\right)$ with $n
\geqslant 1$, ${\mathcal F}'={\mathcal F}^\xi:=\{a \in {\mathcal F}\,|\,\,{\rm
supp\,}a \subset [T-\xi,T]\}$ with a fixed positive $\xi<T$. Also,
let $e^\xi: {\mathcal F}^\xi \to {\mathcal F}$ be the embedding, and
${{\mathcal P}}={{\mathcal P}}^\xi$ the projector specified in Lemma
\ref{lemma cPxi via blocks}. In addition, assume that the
isomorphism $C$ in the Lemma is of the form
 \begin{equation}\label{general C}
(Ca)(t)=\gamma a(t)+\int_0^Tc(t,s) a(s)\,ds\,, \qquad 0\leqslant t
\leqslant T
 \end{equation}
with a constant $\gamma>0$ and the matrix kernel $c(t,s)$ \,
continuous outside the dia\-gonal $t=s$. Fix an $a\in {\mathcal F}$;
let $C^\xi:=(e^\xi)^*C e^\xi$ be the block of $C$ in ${\mathcal
F}^\xi$. By (\ref{repres cP via blocks}), the projection ${{\mathcal
P}}^\xi a=:a^\xi$ satisfies
 $
C^\xi (e^\xi)^*a^\xi\,=\,(e^\xi)^*C a
 $
that is now equivalent to the equation
 \begin{equation}\label{int eqn a xi}
\gamma a^\xi(t)+\int_{T-\xi}^Tc(t,s) a^\xi(s)\,ds=\gamma
a(t)+\int_0^Tc(t,s) a(s)\,ds\,, \quad T-\xi \leqslant t \leqslant
T\,.
 \end{equation}
By the general integral equation theory, the inverse to $C^\xi$ in
${\mathcal F}^\xi$ is of the form
 \begin{equation}\label{inv C xi}
\left([C^\xi]^{-1}b\right)(t)=\gamma^{-1}b(t)+\int_{T-\xi}^T
k^\xi(t,s) b(s)\,ds\,, \qquad T-\xi \leqslant t \leqslant T
 \end{equation}
with the kernel $k^\xi$ continuous outside the diagonal $t=s$ and
continuously depending on $\xi$. For the kernels, the relation
$[C^\xi]^{-1}C^\xi=\mathbb I$ is equivalent to the equation
 \begin{equation}\label{int eqn for kernels}
\gamma^{-1} c(t,s)+\gamma k^\xi(t,s)+\int_{T-\xi}^T k^\xi(t,\eta)
c(\eta,s)\,d\eta=0\,, \quad T-\xi \leqslant s,t \leqslant T\,,
 \end{equation}
which determines $k^\xi$ via $c$. Also, if $c$ is $C^1$-smooth
outside the diagonal $t=s$ then the solution $k^\xi$ is of the
same smoothness character.

 \begin{lemma}\label{lemma cPxi integral repres}
The representation
 \begin{equation}\label{repres integral Pxi}
\left({{\mathcal P}^\xi} a\right)(t)=\begin{cases} 0, & 0 \leqslant t
< T-\xi\\a(t)+\int_0^{T-\xi}p^\xi(t,s) a(s)\,ds, & T-\xi \leqslant
t\leqslant T
                        \end{cases}
 \end{equation}
holds with a continuous matrix kernel $p^\xi$.
 \end{lemma}
$\square$\quad The equality in the first line of (\ref{repres
integral Pxi}) is evident. Let us justify the second line.

The embedding $e^\xi$ extends ${\mathbb C}^n$-valued functions
from $[T-\xi, T]$ to $[0,T-\xi]$ by zero, whereas $(e^\xi)^*$
restricts functions on the segment $[T-\xi, T]$.

Applying $[C^\xi]^{-1}$ to (\ref{int eqn a xi}), one has
 \begin{align}
\notag & a^\xi(t)\overset{(\ref{inv C
xi})}=a(t)+\\
\notag & +\int_0^T\gamma^{-1}c(t,s) a(s)\,ds +\int_{T-\xi}^T
k^\xi(t,\eta)\left[\gamma
a(\eta)+\int_0^T c(\eta,s)a(s)\,ds\right]\,d\eta=\\
\notag
 &=a(t)+\left[\int_0^{T-\xi}\!+\!\int_{T-\xi}^T\right]\gamma^{-1}c(t,s)a(s)\,ds+\int_{T-\xi}^T\gamma
k^\xi(t,s)a(s)\,ds+\\
\notag
&+\left[\int_0^{T-\xi}\!+\!\int_{T-\xi}^T\right]\left\{\int_{T-\xi}^Tk^\xi(t,\eta)c(\eta,s)\,d\eta\right\}a(s)\,ds=\\
\notag & = a(t)+
\int_0^{T-\xi}\left[\gamma^{-1}c(t,s)+\int_{T-\xi}^T k^\xi(t,\eta)
c(\eta,s)\,d\eta\right]a(s)\,ds+\\
\notag & +\int_{T-\xi}^T\left\{\gamma^{-1}c(t,s)+\gamma
k^\xi(t,s)+\int_{T-\xi}^Tk^\xi(t,\eta)c(\eta,s)\,d\eta\right\}a(s)\,ds=\\
\notag & = a(t)+\int_0^{T-\xi}p^\xi(t,s) a(s)\,ds\,, \qquad T-\xi
\leqslant t \leqslant T
 \end{align}
since the expression in the braces vanishes by (\ref{int eqn for
kernels}). The kernel
 \begin{equation}\label{VV}
p^\xi(t,s):=\gamma^{-1}c(t,s)+\int_{T-\xi}^T k^\xi(t,\eta)
c(\eta,s)\,d\eta
 \end{equation}
is continuous since the point $(t,s)$ doesn't reach the diagonal
as $0<s<T-\xi$ and $T-\xi<t<T-\xi$. \quad $\blacksquare$
 \begin{corollary}
The kernel $p^\xi$ of the projector ${{\mathcal P}}^\xi$ satisfies
  \begin{equation}\label{p(t,t)}
p^\xi(T-\xi, T-\xi)=- \gamma k^\xi(T-\xi, T-\xi)\,.
  \end{equation}
 \end{corollary}
To establish (\ref{p(t,t)}) it suffices to put $t=s=T-\xi$ in
(\ref{VV}) and use (\ref{int eqn for kernels}).

\subsubsection*{Operator $W$}
Regarding $\xi \in (0,T]$ as a parameter, we have a family of
projectors ${{\mathcal P}}^\xi$. The family determines an operator $W$, which
is the key object of our constructions.

Let ${\mathcal H}:=L_2(\Omega^T; {\mathbb C}^2)$ be the space of
functions of variable $x:\,0 \leqslant x \leqslant T$, $\sigma$ a
matrix obeying
 \begin{equation}\label{sigma*sigma=gamma}
\sigma^* \sigma\,=\,
 \begin{pmatrix}
\gamma & 0\\0 & \gamma
 \end{pmatrix}\,,
 \end{equation}
where $(\sigma^*)_{ij}=\overline{\sigma_{ji}}$. Define an operator
$W: {\mathcal F} \to {\mathcal H}$,
 \begin{align}
\notag & \left(W a\right)(x):=\sigma \left({{\mathcal P}}^x
a\right)\big|_{t=T-x-0}=\\
\notag & \overset{(\ref{repres integral
Pxi})}=\sigma a(T-x)+\int_0^{T-x}\sigma p^x(T-x,s) a(s)\,ds=\\
 \label{def W} &=\sigma a(T-x)+\int_x^T\hat w(x,s)
 a(T-s)\,ds\,, \qquad x \in \Omega^T\,,
 \end{align}
with the continuous matrix kernel
 \begin{equation}\label{hat w}
\hat w(x,s)\,:=\,\sigma p^x(T-x, T-s)\,, \qquad 0\leqslant x
\leqslant s \leqslant T\,.
 \end{equation}
satisfying
 \begin{equation}\label{diag hat w}
\hat w(x,x)\overset{(\ref{p(t,t)})}=\, - \gamma\sigma \, k^x(T-x,
T-x)\,, \qquad x \in \Omega^T
 \end{equation}
and
\begin{equation}\label{hat w(0,s)}
\hat w(0,s)\overset{(\ref{hat w})}=\,\sigma
p^0(T,T-s)\overset{(\ref{VV})}=\gamma^{-1}\sigma \,c(T,T-s)\,,
\qquad 0\leqslant s \leqslant T\,.
 \end{equation}
{\bf Remark 1.}\quad A simple analysis shows: if kernel $c(t,s)$
in (\ref{general C}) is $C^1$-smooth outside the diagonal $t=s$
then kernel $\hat w$ is $C^1$-smooth in the domain of its
definition. As a result, $\hat w(x,x)$ is $C^1$-smooth for $x \in
\Omega^T$.

Let $P^\xi$ be the orthogonal projector in ${\mathcal H}$ onto
${{\mathcal H}^\xi}:=\{y \in {\mathcal H}\,|\,\,{\rm supp\,}y \subset
\Omega^\xi\}$. It cuts off functions on the segment $\Omega^\xi$:
 \begin{equation*}
\left(P^\xi y\right)(x)=
    \begin{cases}
y(x), & 0\leqslant x \leqslant \xi\\0, & \xi< x \leqslant T
    \end{cases}
 \end{equation*}
and satisfies $(P^\xi)^*=P^\xi$. Also, introduce $Y: {\mathcal
F}\to{\mathcal H}$,
 \begin{equation*}
(Ya)(x)\,:=\,a(T-x)\,, \qquad x \in \Omega^T\,,
 \end{equation*}
so that $YY^*={\mathbb I}_{\mathcal H}$. The following is the basic
properties of operator $W$.
\begin{lemma}\label{lemma W cP = P W}
The relations
\begin{equation}\label{W cP = P W}
W {{\mathcal P}}^\xi\,=\,P^\xi W\,, \qquad C=W^*W
 \end{equation}
hold for $0<\xi \leqslant T$.
 \end{lemma}
$\square$ \quad
By the representation in (\ref{def W}), $W$ is of
the form `isomorphism + Volterra type operator'. Therefore, $W$ is
an isomorphism from ${\mathcal F}$ onto ${\mathcal H}$. Also, the same
representation easily implies $W{\mathcal F}^\xi = {\mathcal H}^\xi$. So,
in terms and notation of Lemma \ref{lemma cPxi via blocks}, $W:
{\mathcal F}_C \to {\mathcal H}$ is a unitary operator, which maps ${\mathcal
F}^\xi$ onto ${\mathcal H}^\xi$. Such an operator interlaces the
corresponding projectors: $W {{\mathcal P}}^\xi=P^\xi W$.

The operator $A:= WC^{-1}W^*$ is an isomorphism in ${\mathcal H}$. It
satisfies
 \begin{align}
\notag & A P^\xi=WC^{-1}W^*P^\xi=WC^{-1}(P^\xi
W)^*=WC^{-1}({{\mathcal P}}^\xi)^*W^*\overset{(\ref{repres cP via
blocks})}=\\
\notag & = W {{\mathcal P}}^\xi C^{-1}W^*=P^\xi W C^{-1}W^*=P^\xi A\,,
\end{align}
so that $A P^\xi=P^\xi A$ for all $0<\xi {\leqslant} T$. By the use of
(\ref{general C}) and (\ref{def W}), one derives
 $$
A=\left[\sigma Y+B\right]\left[\gamma^{-1}{\mathbb I}_{\mathcal
F}+K\right]\left[ Y^* \sigma^* +B^*\right]=\sigma
\gamma^{-1}\sigma^*\, {\mathbb I}_{\mathcal
H}+D\,\overset{(\ref{sigma*sigma=gamma})}=\,{\mathbb I}_{\mathcal
H}+D,
 $$
where $B,\,K,\,D$ are the {\it compact} compact operators, $D$
commuting with {\it all} $P^\xi$. Since ${\rm dim\,}{\mathcal
H}=\infty$, the latter is possible only if $D=\mathbb O$. Thus,
$A={\mathbb I}_{\mathcal H}$ that is equivalent to the second relation
in (\ref{W cP = P W}). \quad $\blacksquare$
\smallskip

\noindent{\bf Remark 2.}\quad As one can easily check, the
operator $F:=Y^*W=\sigma {\mathbb I}_{\mathcal F}+ $ $\tilde B$
satisfies $F^*F=C$ and $F{\mathcal F}^\xi={\mathcal F}^\xi$ for all $\xi$.
This means that $F$ provides a {\it triangular factorization} to
$C$ in ${\mathcal F}$ with respect to the family of subspaces $\{{\mathcal
F}^\xi\}$ \cite{GK}, \cite{Dav}. Such a factorization is unique up
to the choice of the matrix $\sigma$ satisfying
(\ref{sigma*sigma=gamma}). Also, written in the form $W{{\mathcal
P}}^\xi W^{-1}=P^\xi$ the first relation in (\ref{W cP = P W})
shows that the map $W: {\mathcal F}_C \to {\mathcal H}$ {\it diagonalizes}
the projector family $\{{{\mathcal P}}^\xi\}$ in the sense of the
Spectral Theorem \cite{BSol}.

\subsection{Solving IP}\label{subsec Solving IP}

\subsubsection*{Determination of potential}
Applying the above-developed theory to system ${{{\beta}^T}}$, we
take ${\mathcal F}={{{\mathcal F}^T_{\beta}}}$, ${\mathcal F}^\xi={{{\mathcal
F}^{T, \xi}_{\beta}}}$, ${\mathcal H}={{\mathcal H}^T}$, ${\mathcal
H}^\xi\equiv {\mathcal H}^\xi$, and $C=C^T$  (see (\ref{CT beta via
R2T alpha}) and (\ref{general C}) with $\gamma=2$). Further,
${{\mathcal P}}^\xi$ projects in ${{{\mathcal F}^T_{\beta}}}$ onto
${{{\mathcal F}^{T, \xi}_{\beta}}}$ in parallel to $(C^T)^{-1}[{{{\mathcal
F}^T_{\beta}}} \ominus {{{\mathcal F}^{T, \xi}_{\beta}}}]$, projects
in $({{{\mathcal F}^T_{\beta}}})_{C^T}$ onto ${{{\mathcal F}^{T,
\xi}_{\beta}}}$ orthogonally, and is represented by (\ref{repres
integral Pxi}) via the kernel (\ref{VV}) satisfying
(\ref{p(t,t)}).

By (\ref{U^xi=H^xi in beta}) \footnote{i.e., owing to the
controllability of ${{{\beta}^T}}$!}, one has $W^T{{\mathcal P}}^\xi=P^\xi W^T$. For a
continuous $a \in {{{\mathcal F}^T_{\beta}}}$, the latter implies
 \begin{align*}
& \left(W^Ta\right)(\xi)=\left(P^\xi W^T
a\right)(T-\xi-0)=\left(W^T {{\mathcal P}}^\xi
a\right)(T-\xi-0)\overset{(\ref{WT repres})}=\\
& =\varkappa\left({{\mathcal P}}^\xi a\right)(T-\xi-0)+\int_\xi^T \check
w(\xi,s) \left({{\mathcal P}}^\xi a\right)(T-s)\,ds=\\
&=\varkappa\left({{\mathcal P}}^\xi
a\right)(T-\xi-0)\overset{(\ref{def W})}=\left(W a\right)(\xi)\,,
\qquad 0<\xi \leqslant T\,,
 \end{align*}
where $W$ corresponds to $\sigma=\varkappa$ (see (\ref{def
varkappa and check w}), (\ref{WT repres})), the integral vanishing
by virtue of ${\rm supp\,}{{\mathcal P}}^\xi a \subset [T-\xi,T]$.
Thus, in our case, $W^T$ is identical to $W$. Hence, their kernels
$\check w$ and $\hat w$ coincide. The latter implies
 \begin{equation}\label{diag check w}
\check w(x,x)\overset{(\ref{diag hat w})}=\, - \gamma\varkappa \,
k^x(T-x, T-x)\,, \qquad x \in \Omega^T\,.
 \end{equation}

To solve the IP it suffices just to summarize these considerations
in the form of the procedure.
\smallskip

\noindent{\bf$1.$}\,\,\,Given the response function
$r_\alpha(t),\,\,0\leqslant t \leqslant 2T$ of the system
${{{\alpha}^T}}$, determine the matrix-kernel $c^T$ by
(\ref{kernel c in detail in forward problem}).
\smallskip

\noindent{\bf$2.$}\,\,\,Solve the family of the {\it linear
integral equations} (\ref{int eqn for kernels}) for $c=c^T$ and
get the kernels $k^\xi,\,\,\,0<\xi\leqslant T$. The solvability is
guaranteed by the positive-definiteness of $C^T$ (and all of its
blocks $C^\xi$). By standard integral equations theory arguments,
the solution $k^\xi$ is of the same smoothness as $c^T$, i.e., is
$C^1$-smooth outside the diagonal $t=s$.
\smallskip

\noindent{\bf$3.$}\,\,\,find the matrix $\check w(x,x)$ by
(\ref{diag check w}), take its first column and recover the
potential by (\ref{p+iq=w2-iw1}): $\Im w_1(x,x)+ \Re
w_2(x,x)=p(x),\,\,-\Re w_1(x,x)+ \Im w_2(x,x)=q(x)$ for $x \in
\Omega^T$.
\smallskip

In this procedure, the equations (\ref{int eqn for kernels}) play
the role of relevant version of the classical
Gelfand-Levitan-Krein-Marchenko equations \cite{BM_1} Also, note
that to solve (\ref{int eqn for kernels}) is, in fact, to invert
the connecting operator $C^T$. So, to solve the IP is to find
$[C^T]^{-1}$ that is a `paradigm' of the BC-method.
\smallskip

\noindent{\bf Remark 3.}\,\,The determination $r \Rightarrow V$ by
the procedure 1.-- 3. is quite transparent and convenient for
stability analysis. One can easily show that a Lipschitz-type
stability does hold: small $C^1$-variations of $r$ imply small
$C^1$-variations of $V$, the corresponding Lipschitz constant
being determined by the (positive) low bound of the spectrum of
$C^T$.

Apart from determination of coefficients (potential), the
BC-method provides one more option that is referred to as {\it
visualization of waves} \cite{BIP07}, \cite{BM_1}. Namely,
possessing the inverse data, one can construct the projectors
${{\mathcal P}}^\xi$ and then recover the control operator by $(W^T
a)(x)\overset{(\ref{def W})}=\varkappa \left({{\mathcal P}}^x
a\right)\big|_{t=T-x-0}$, i.e., determine the solutions (states)
without solving the forward problem.

\subsubsection*{Characterization of data} It remains to prove
Theorem \ref{Th Charact}, i.e., provide the solvability conditions
for the IP under consideration.
\smallskip

\noindent {\bf Necessity.}\quad If the function $r$ is a response
function of a system ${{{\alpha}^T}}$ (so that $r=r_\alpha$) then the
operator $C$ given by (\ref{cal C repres}) coincides with the
connecting operator $C^T$: see (\ref{CT beta via R2T alpha}). The
latter is a positive definite isomorphism in ${{{\mathcal F}^T_{\beta}}}$.
\smallskip

\noindent {\bf Sufficiency.}\quad {\bf $1$.}\,\,\,Let $r$ be such
that the operator $C$ defined by (\ref{cal C repres}) is
po\-sitive definite in ${\mathcal F}:=L_2([0,2T];{\mathbb C}^2)$.
Written in detail with regard to (\ref{cal R repres}) and
(\ref{cal RT repres}), it takes the form
 \begin{equation*}
\left(Ca\right)(t)\,=\,\gamma a(t)+\int_0^T c(t,s) a(s)\,ds\,,
\qquad 0\leqslant t \leqslant T
 \end{equation*}
with $\gamma =2$ and the matrix kernel elements
 \begin{align}
\notag & c_{11}(t,s)(t)=-i\,[r(t-s)-\bar r(s-t)]\,, \quad
c_{12}(t,s)=-i\,\bar r(2T-t-s)\,,\\
& c_{21}(t,s)= i\,r(2T-t-s)\,,\quad c_{22}(t,s)(t)=i\,[\bar
r(t-s)-r(s-t)]\label{kernel c in detail}
 \end{align}
(recall Convention \ref{Conv 1}). Note the relation
 \begin{equation*}
c(T,T-s)=
   \begin{pmatrix}
-i r(s) & -i \bar r(s)\\
i r(s) & i \bar r(s)
   \end{pmatrix}\,, \qquad 0\leqslant s \leqslant T
 \end{equation*}
following from (\ref{kernel c in detail}), and its consequence
\begin{equation}\label{gamma{-1}kappa c(T,T-s)}
\gamma^{-1}\varkappa\, c(T,T-s) =
   \begin{pmatrix}
0 & 0 \\
r(s) & \bar r(s)
   \end{pmatrix}\,, \qquad 0\leqslant s \leqslant T\,.
 \end{equation}
\smallskip

{\bf $2$.}\,\,\,In ${\mathcal F}$ introduce the classes ${{\mathcal M}}^T:=\{a
\in C^1\left([0,T]; {\mathbb C}^2\right)\,|\,\,a(0)=a'(0)=0\}$ and
 $$
{{\mathcal M}}^T_0:= \left\{a=\begin{pmatrix}a_1\\a_2 \end{pmatrix}\in
{{\mathcal M}}^T\,\bigg|\,\,a_1(T)+a_2(T)=0\right\}\,.
 $$
Recall that $Q:=\begin{pmatrix}1&0\\0&-1\end{pmatrix}$. A simple
analysis, which uses nothing but integration by parts and a
con\-vo\-lution form (\ref{kernel c in detail}) of $c_{kl}$, leads
to the relation
 \begin{equation}\label{symmetry -iQd/dt}
\left(C\left[-iQ\frac{da}{dt}\right], b\right)_{\mathcal
F}\,=\,\left(Ca,\left[-iQ\frac{db}{dt}\right]\right)_{\mathcal F}
\qquad \text{for} \,\,\,a,b \in {{\mathcal M}}^T_0\,.
 \end{equation}

Recall that the space ${\mathcal F}_C$ is defined in sec \ref{subsec A
portion of Operator Theory}. Introduce the operator $D_0: {\mathcal
F}^T \to {\mathcal F}^T,\,\,\,{\rm Dom\,} D_0={{\mathcal M}}^T_0,\,\,D_0:=-i
Q\frac{d}{dt}$. Then (\ref{symmetry -iQd/dt}) takes the form \begin{equation}
\label{symmetry D 0} \left(D_0 a, b\right)_{{\mathcal
F}_C}\,=\,\left(a, D_0 b\right)_{{\mathcal F}_C}\,,
\end{equation}
and shows that
$D_0$ is a symmetric densely defined operator in ${\mathcal F}_C$.
\smallskip

{\bf $3$.}\,\,\, In accordance with the scheme of sec \ref{subsec
A portion of Operator Theory}, operator $C$ determines the
operator $W:{\mathcal F} \to {{\mathcal H}^T}:=L_2(\Omega^T;{\mathbb
C}^2)$, which acts by the rule
 \begin{equation}\label{char def W-1}
\left(W a\right)(x)\overset{(\ref{def W})}=\varkappa
a(T-x)+\int_x^T\hat w(x,s)
 a(T-s)\,ds\,, \qquad x \in \Omega^T\,,
 \end{equation}
and satisfies
 \begin{equation}\label{char C=W*W}
C\,\overset{(\ref{W cP = P W})}=\,W^*W\,.
 \end{equation}
By (\ref{hat w(0,s)}) (with $\sigma=\varkappa$) and
(\ref{gamma{-1}kappa c(T,T-s)}), the kernel $\hat w$ obeys
\begin{equation}\label{ww}
\hat w(0,s)=
\begin{pmatrix}
0 & 0 \\
r(s) & \bar r(s)
   \end{pmatrix}\,, \qquad 0\leqslant s \leqslant T\,.
 \end{equation}

In ${{\mathcal H}^T}$ introduce the classes ${{\mathcal N}}^T:=\{y \in
C^1\left([0,T]; {\mathbb C}^2\right)\,|\,\,y(T)=y'(T)=0\}$ and
 $$
{{\mathcal N}}^T_0:= \left\{y=\begin{pmatrix}y_1\\y_2 \end{pmatrix}\in
{{\mathcal N}}^T\,\bigg|\,\,y_1(0)=0\right\}\,.
 $$
The following relations are valid:
\begin{equation}\label{W M0=N0}
W {{\mathcal M}}^T\,=\,{{\mathcal N}}^T\,, \quad W {{\mathcal
M}}^T_0\,=\,{{\mathcal N}}^T_0\,.
 \end{equation}
Indeed, by $C^1$-smoothness of kernel $\hat w$, operator $W$ maps
$C^1([0,T]; {\mathbb C}^2)$ onto $C^1(\Omega^T; {\mathbb C}^2)$
isomorphically. As it follows from the form of the integral term
in (\ref{char def W-1}), the relations $a(0)=a'(0)=0$ and
$(Wa)(T)=(Wa)'(T)=0$ are equivalent. At last, $(Wa)_1(0)=y_1(0)$
is equivalent to $a_1(T)+a_2(T)=y_1(0)$ by (\ref{char def W-1})
and (\ref{ww}). Hence, the second equality in (\ref{W M0=N0}) is
also valid.
\smallskip

{\bf $4$.}\,\,\,One can write $W$ in the form
 \begin{equation}\label{W via M}
W\,=\,\varkappa\left[{\mathbb I}+M\right]Y, \quad
(My)(x):=\int_x^Tm(x,s) y(s)\,ds
 \end{equation}
with $Y: {\mathcal F} \to {{\mathcal H}^T},\,\,(Ya)(x)=a(T-x)$ and an
integral Voterra operator $M: {{\mathcal H}^T} \to {{\mathcal H}^T}$ with
the kernel $m:=\varkappa^{-1}\hat w$. Respectively, its inverse
takes the form
 \begin{equation}\label{W^-1 via N}
W^{-1}=Y^{-1}\left[{\mathbb I}+N\right]\varkappa^{-1}, \quad
(Ny)(x):=\int_x^T n(x,s) y(s)\,ds\,,
 \end{equation}
where $N: {{\mathcal H}^T} \to {{\mathcal H}^T}$ is an integral Volterra
operator. Both of the kernels are $C^1$-smooth in $\Delta^T$. The
equality $W W^{-1}=\mathbb I$ is equivalent to
\begin{equation}\label{m+n+int mn=0}
m(t,s)+n(t,s)+\int_s^t m(t,\tau) n(\tau,s)\,d\tau=0\,, \qquad 0
\leqslant s \leqslant t \leqslant T\,.
 \end{equation}
\smallskip

{\bf $5$.}\,\,\ Define an operator $L: {{\mathcal H}^T} \to {{\mathcal
H}^T}$ on ${\rm Dom\,}L={{\mathcal N}}^T$ by
\begin{equation}\label{def L}
L\,:=\,W \left[-iQ\frac{d}{dt}\right] W^{-1}
 \end{equation}
\begin{lemma}\label{lemma W[d/dt]W{-1}}
The representation
 \begin{equation*}
L\,=\,J\frac{d}{dx} + V\,, \qquad V=
  \begin{pmatrix}p &q\\q & -p\end{pmatrix}
 \end{equation*}
holds with the real $p, q \in C^1(\Omega^T)$.
 \end{lemma}
$\square$\quad For a $y \in {{\mathcal N}}^T$, with regard to (\ref{W via M})
and (\ref{W^-1 via N}), one has
 $$
\left[\varkappa^{-1}L\varkappa\right]y=i\left[{\mathbb I}+\int_x^T
m\right]Y\left[-iQ\frac{d}{dt}\right]Y^{-1}\left[{\mathbb
I}+\int_x^T n\right]y\,.
 $$
Applying, step by step, the operators, and integrating by parts
with regard to
 $$
y(T)=0,\quad \frac{d}{dt}Y^{-1}=-Y^{-1}\frac{d}{dx},\quad
YQY^{-1}=Q\,,
 $$
one gets
 \begin{align*}
&
\left(\left[\varkappa^{-1}L\varkappa\right]y\right)(x)=\\
& =iQy'(x)-i\left[m(x,x)Q+Qn(x,x)\right]y(x)+
\int_x^Tj(x,s)y(s)\,ds\overset{(\ref{m+n+int mn=0})}=\\
& = iQy'(x)+i\left[Qm(x,x)-m(x,x)Q\right]y(x)+
\int_x^Tj(x,s)y(s)\,ds\overset{(\ref{m+n+int mn=0})}=\\
& = iQy'(x)+i\, 2\begin{pmatrix}0&m_{12}(x,x)\\-m_{21}(x,x) &
0\end{pmatrix}y(x)+ \int_x^Tj(x,s)y(s)\,ds\,,
 \end{align*}
where $m_{kl}$ are the matrix elements of the kernel $m$. Then,
with regard to $i\varkappa Q \varkappa^{-1}=J$, one easily arrives
at
 \begin{equation}\label{L with int yet}
(Ly)(x)=Jy'(x)+\begin{pmatrix}p(x) &q(x)\\q(x) &
-p(x)\end{pmatrix}y(x)+\int _x^T h(x,s) y(s)\,ds
 \end{equation}
with the continuously differentiable
 $$
p(x)=i[m_{12}(x,x)-m_{21}(x,x)]\,, \quad
q(x)=-[m_{12}(x,x)+m_{21}(x,x)]\,,
 $$
and a continuous kernel $h$.
\smallskip

Show that $p$ and $q$ are real, and $h \equiv 0$. Let
$L_0:=L\big|_{{{\mathcal N}}^T_0}$. The definition (\ref{def L}) and second
relation in (\ref{W M0=N0}) imply $L_0=W D_0 W^{-1}$, whereas $W$
is a unitary operator as a map from ${\mathcal F}_C$ to ${{\mathcal
H}^T}$. By the latter, since $D_0$ is a symmetric densely defined
operator in ${\mathcal F}_C$ (see (\ref{symmetry D 0})), the same
holds for $L_0$ in ${{\mathcal H}^T}$ . Looking at its form (\ref{L
with int yet}), it is easy to recognize that $L_0$ can be
symmetric if and only if $V=V^*$ holds (that is $p=\bar p,\,
q=\bar q$) and the {\it triangular} kernel $h$ vanishes
identically. \quad $\blacksquare$
\smallskip

{\bf $6$.}\,\,\,By this, we have the real symmetric zero-trace
potential $V$ {\it determined by the input function} $r$. The
potential, in turn, determines the system $\beta^T$ with the
control operator $W^T$ characterized by relation (\ref{Wd/dt=LW}).
Comparing the latter with (\ref{def L}), we see that $W^T=W$. For
the connecting operator of $\beta^T$, one has
 $$
C^T:=(W^T)^*W^T=W^*W\overset{(\ref{char C=W*W})}= C
 $$
that, in particular, implies $c^T_{21}(t,s)=c_{21}(t,s)$.
Comparing (\ref{kernel c in detail in forward problem}) with
(\ref{kernel c in detail}) and taking there $t=s=T-\tau$, we
arrive at $r_\alpha(\tau)=r(\tau),\,\,0\leqslant \tau \leqslant
2T$.
\smallskip

Thus, there exists the system ${{{\alpha}^T}}$, which possesses
the response function $r_\alpha=r$. Theorem \ref{Th Charact} is
proved.

\noindent{\bf Key words:} one-dimensional dynamical Dirac system,
controllability, determination of potential,  characterization of
inverse data.

\noindent{\bf MSC:}\,\,\,35R30, 35Bxx, 35Lxx, 35Qxx.

\end{document}